\documentclass{article}

\usepackage{CJK}
\usepackage{CJK,CJKnumb,CJKulem,times,dsfont,ifthen,mathrsfs,latexsym,amsfonts,color}
\usepackage{amsmath,amsthm,makeidx,fontenc,amssymb,bm,graphicx,psfrag,listings,curves,extarrows}

\usepackage[utf8]{inputenc}
\usepackage[english]{babel}
\usepackage{multicol,graphicx,color}
\usepackage{pslatex}
\usepackage{authblk}
\usepackage{amsthm}
\usepackage{cite}
\usepackage{amsmath}
\usepackage{amssymb}
\usepackage{latexsym}
\usepackage{lscape}
\usepackage{epsfig}
\usepackage{amsfonts}
\usepackage{mathrsfs,bbm}
\usepackage[
hmarginratio={1:1},     
vmarginratio={1:1},     
textwidth=15cm,        
textheight=21cm,
heightrounded,          
]{geometry}

\usepackage{tikz}

\usepackage[
linktocpage=true,colorlinks,citecolor=magenta,linkcolor=blue,urlcolor=magenta]{hyperref}
\makeatletter

\usepackage[showonlyrefs]{mathtools}  
\mathtoolsset{showonlyrefs=true}

\newtheorem{theorem}{Theorem}[section]

\newtheorem{definition}[theorem]{Definition}

\newtheorem{lemma}[theorem]{Lemma}
\newtheorem{proposition}[theorem]{Proposition}
\newtheorem{remark}[theorem]{Remark}

\theoremstyle{definition} \theoremstyle{remark}
\numberwithin{equation}{section}

\newcommand{\R}{\mathbb{R}}
\newcommand{\la}{\lambda}
\newcommand{\ep}{\epsilon}
\newcommand{\C}{\mathcal{C}}
\newcommand{\s}{\mathbb{S}^{N-1}_+}
\newcommand{\cs}{\mathbb{S}^{N-1}}
\newcommand{\A}{\mathcal{A}}
\newcommand{\cE}{\mathscr{E}}
\newcommand{\cX}{\mathcal X}
\newcommand{\cY}{\mathcal Y}
\newcommand{\cF}{\mathcal F}

\begin{document}
\title{\bf Boundary singularities for weak solutions of semilinear elliptic problems throughout the subcritical range
}

	\author{{\bf Linjie Song}\footnote{songlinjie18@mails.ucas.edu.cn, linjie.song@univ-fcomte.fr} }

	\affil{{\small Universit\'e Marie et Louis Pasteur, CNRS, LmB (UMR 6623), Besan\c{c}on F-25000, France}}
	
	\author{{\bf Wenming Zou}\footnote{zou-wm@mail.tsinghua.edu.cn} }
	
	\affil{{\small Department of Mathematical Sciences, Tsinghua University, Beijing 100084, China}}

	\date{}	
	\maketitle
	
	\begin{abstract}

		Let $\Omega \subset \R^N$ ($N \ge 2$) be a bounded domain with smooth boundary $\partial \Omega$. We construct positive very weak solutions for the problem $\Delta u + u^p = 0$ in $\Omega$ for all $p \in [(N+1)/(N-1),(N+2)/(N-2))$ when $N \ge 3$ and for all $p \ge (N+1)/(N-1)$ when $N = 2$. These solutions vanish in a suitable trace sense on $\partial \Omega$ but are singular at prescribed isolated points. This gives a completely affirmative answer to a conjecture of del Pino-Musso-Pacard.

\noindent

		\vskip0.1in
		{\small \noindent \text{\bf Keywords:} Prescribed boundary singularities; Very weak solution; Lane-Emden equation
		\vskip0.1in
		\noindent\text{\bf Mathematics Subject Classification:} 35A01, 35B09 
		\vskip0.1in
		\noindent \textbf{Statements and Declarations:} The authors have no relevant financial or non-financial interests to disclose.
		\vskip0.1in
		\noindent \textbf{Data availability:} Data sharing is not applicable to this article as no datasets were generated or analysed during the current study.	
     	}	
	\end{abstract}
	

	\bigskip
	
	\section{Introduction}
	
	Let $\Omega \subset \R^N$ ($N \ge 2$) be a bounded domain with smooth boundary $\partial \Omega$. We study the following Lane-Emden equation
	\begin{align} \label{eqmain}
		\left\{
		\begin{aligned}
			& -\Delta u = u^p \quad & \text{ in } \Omega, ~~~\\
			& u > 0 \quad & \text{ in } \Omega, ~~~\\
			& u = 0 \quad & \text{ on } \partial \Omega,
		\end{aligned}
    	\right.
	\end{align}
    where $p > 1$, and we are interested in finding {\it very weak solutions} within the meaning given in \cite{DMP}, namely
	\[
	u \in L^1(\Omega), \quad \text{dist}(x,\partial \Omega)u^p \in L^1(\Omega),
	\]
	and 
	\[
	\int_\Omega (u\Delta v + u^p v)dx = 0 \quad \text{for all } v \in C^2(\overline{\Omega}) \text{ with } v = 0 \text{ on } \partial \Omega.
	\]
    
    Let 
    \[
     p_c = \frac{N+1}{N-1}; \qquad p_s = \frac{N+2}{N-2}, ~ N \ge 3, ~~\text{ and }~~ p_s = +\infty, ~ N = 2.
    \]
    When $1 < p < p_c$, it follows from the results in \cite{BT,QS} that a very weak solution is actually in $H_0^1(\Omega)$ and it is a weak solution in the usual sense. As far as boundary regularity of very weak solutions is concerned, results in \cite{Soup} imply that the exponent $p = p_c$ is critical. For exponents slightly larger than or equal to $p_c$, del Pino, Musso and Pacard \cite{DMP} constructed solutions to problem \eqref{eqmain} with prescribed singularities on the boundary. More precisely, they introduced the following definition:    
    \begin{definition}
    	Let $u(x)$ be a function defined in $\Omega$ and $y \in \partial \Omega$. We say that 
    	\begin{center}
    		$u(x) \to \ell$ ~~ as $x \to y$ nontangentially
    	\end{center}
        if
        \begin{align*}
        	\lim_{\Gamma_\alpha(y) \ni x \to y} u(x) = \ell \quad \text{for all } \alpha \in \left[0,\frac{\pi}{2}\right),
        \end{align*}
        where $\Gamma_\alpha(y)$ denotes the cone with vertex $y$, and angle $\alpha$ with respect to its axis, the inner normal to $\partial \Omega$ at $y$.
    \end{definition}
    \noindent Their main result (\cite[Theorem 1.1]{DMP}) reads:
    \begin{theorem} \label{thm:dmp}
    	There exists a number $p_N > p_c$ such that, given $p \in [p_c,p_N)$ and given points $y_1, y_2, \cdots, y_k \in \partial \Omega$, there exist very weak solutions $u$ to problem \eqref{eqmain} such that $u \in C^2(\overline{\Omega} \backslash \{y_1, y_2, \cdots, y_k\})$ and
    	\begin{center}
    		$u(x) \to +\infty$ ~~ as $x \to y_i$ nontangentially, for all $i = 1,\cdots, k$.
    	\end{center}
    \end{theorem}   
   
    Furthermore, del Pino, Musso and Pacard \cite{DMP} conjectured that \textit{the above result should hold for $p \in [p_c,p_s)$}. In this paper, we give a completely positive answer to this conjecture. The result is as follows.
    \begin{theorem} \label{thmmain}
    	Given $p \in [p_c,p_s)$ and given points $y_1, y_2, \cdots, y_k \in \partial \Omega$, there exist very weak solutions $u$ to problem \eqref{eqmain} such that $u \in C^2(\overline{\Omega} \backslash \{y_1, y_2, \cdots, y_k\})$ and
    	\begin{center}
    		$u(x) \to +\infty$ ~~ as $x \to y_i$ nontangentially, for all $i = 1,\cdots, k$.
    	\end{center}
    \end{theorem}

	In fact, the very weak solutions shown in above theorem are infinitely many, see details in its proof.

    In the construction of the solutions predicted by Theorem \ref{thmmain}, the basic cells will be positive solutions of an equation defined on the half-space,
    \[
    \R^N_+ := \big\{x = (x_1,\cdots,x_N): x_N > 0\big\},
    \]
    which vanish on its boundary and have a singularity at the origin. For this, we consider the following problem:
	\begin{align} \label{eqhalfspace}
	\left\{
	\begin{aligned}
		& -\Delta u = u^p \quad & \text{ in } \R^N_+ \backslash \{0\}, ~~~\\
		& u > 0 \quad & \text{ in } \R^N_+, ~~~~~~~~~~~\\
		& u = 0 \quad & \text{ on } \partial \R^N_+\backslash \{0\},
	\end{aligned}
	\right.
   \end{align}     
   with $p > p_c$. 
   
   Let us denote by $\mathbb{S}^{N-1}_+$ the half sphere
   \[
   \mathbb{S}^{N-1}_+ := \big\{z = (z_1,\cdots,z_N) \in \mathbb{S}^{N-1}: z_N > 0\big\}.
   \]
   We further set
   \[
   a_p = \frac{4}{p-1} - (N-2), \quad b_p = \frac{2}{p-1}\left(N-2 - \frac{2}{p-1}\right).
   \]
   For $p \in (p_c,p_s)$, we have $a_p > 0$ and $b_p + N-1 > 0$. Let $\phi_p$ denote the unique solution of
   \begin{align} \label{eqhalfsphere}
   	\left\{
   	\begin{aligned}
   		& -\Delta_{\mathbb{S}^{N-1}} \phi + b_p \phi = \phi^p \quad & \text{ in } \mathbb{S}^{N-1}_+, ~~~\\
   		& \phi > 0 \quad & \text{ in } \mathbb{S}^{N-1}_+, ~~~\\
   		& \phi = 0 \quad & \text{ on } \partial \mathbb{S}^{N-1}_+.
   	\end{aligned}
   	\right.
   \end{align}
   Here $\Delta_{\mathbb{S}^{N-1}}$ designates the Laplace-Beltrami operator in $\mathbb{S}_+^{N-1}$. Note that $N-1$ is the first eigenvalue of $-\Delta_{\mathbb{S}^{N-1}}$ on the half sphere under Dirichlet boundary conditions. The corresponding eigenfunction is given by 
   \begin{equation}\label{eq:phi1}
   	\varphi_1(z)=\frac{z_N}{\|z_N\|_{L^2(\s)}}.
   \end{equation}
   See \cite[Theorem 1.1 and Section 3]{BPV} for the existence and uniqueness of $\phi_p$. A function is called {\it axial} if it is invariant under the action of $O(N-1)$ on the first $N-1$ coordinates. As pointed out in \cite[p.192]{BPV}, by an adaptation of the moving planes method to $\mathbb{S}^{N-1}$, it can be shown that $\phi_p$ depends only on the geodesic distance to the north pole
   \[
   \theta = \arccos z_N
   \]
   and is decreasing with respect to $\theta$. Thus $\phi_p$ is an axial function.
   
   Now we can state our result on the half space.   
   \begin{proposition} \label{prophalfspacesolution}
   	Given $p \in (p_c,p_s)$, there exists a solution $u_p$ to problem \eqref{eqhalfspace} such that 
   	\begin{align} \label{eq:originratio}
   		u_p(x) = |x|^{-\frac{2}{p-1}}\phi_p(z_N)\left(1 + o(1)\right) \quad \text{as } x \to 0,
   	\end{align}
    and
    \begin{align} \label{eq:infinityratio}
    	u_p(x) = |x|^{1-N}z_N\left(1 + o(1)\right) \quad \text{as } |x| \to +\infty.
    \end{align}
   	In addition, we have the pointwise estimate
   	\begin{align} \label{eq:globalbound}
   		0 < u_p(x) \le C_p|x|^{-\frac{2}{p-1}}z_N, \quad x \in \R^N_+.
   	\end{align}  
   \end{proposition}

   There are three distinct ingredients in the proof of the half-space existence result. Firstly, the angular equilibrium has exactly one negative axial eigenvalue and no zero axial eigenvalue, see Subsection~\ref{subangular}. The energy identity, together with subcritical blow-up bounds, makes the set of connections compact after a scalar phase condition, see Subsections~\ref{subprioribounds}, \ref{subcompactness}. Finally, the cylindrical linearization has Fredholm index one, and the phase condition lowers this index to zero, see Subsection~\ref{subindex}. Then we complete the proof of Proposition~\ref{prophalfspacesolution} in Subsection~\ref{subproof1}, by studying the zero set of the map given in \eqref{eq:Fmap}. When $p$ is close to $p_c$, the zero with the phase condition is unique and regular. For $p$ away from $p_c$, we use a continuation argument to prove that the zero set is always nonempty. In Appendix~\ref{sec:continuationlemma}, we state and prove the continuation principle, which is a local form of the mod-$2$ Fredholm degree principle of \cite{PR}.

   It is worth to remark that Proposition \ref{prophalfspacesolution} generalizes the result in \cite[Proposition 1.1]{DMP} and this provides a completely affirmative answer to {\it Open problem 1} proposed by del Pino, Musso and Pacard \cite{DMP}. We also remark that the existence result in Proposition \ref{prophalfspacesolution} alone is not enough for a gluing argument shown in \cite{DMP} for all $p \in (p_c,p_s)$ and we will develop a new gluing argument to complete the proof of Theorem \ref{thmmain}.

   For the bounded domain case, the proof starts with a degree argument. In Subsection~\ref{subdegree}, we prove that the compact set of connections with the phase condition has local mod-2 degree one in the full function space, with the axial restrictions removed. Then, in Subsection~\ref{subinner-outer}, we couple rescaled cylindrical profiles near the prescribed boundary points to a regular outer correction. The singular nonlinear terms can be removed from the outer equation due to an exact cancellation, allowing us to solve this equation by contraction. Substituting this correction into the inner equations, we obtain a finite product of the half-space problems with a $C^1$-small perturbation. The degree remains nonzero, yielding the required solution without assuming non-degeneracy of the connections at the prescribed exponent. Finally, the proof of Theorem~\ref{thmmain} is completed in Subsection~\ref{subproof2}.

   After completing the paper, we learned that Ao, Katayama and Wei in a recent preprint \cite{AKW} also solved the conjecture of del Pino-Musso-Pacard in dimensions $N \ge 3$ for $p$ in the interval $(p_c,p_s)$ outside a locally finite set. Their construction is also based on the half-space result stated in Proposition~\ref{prophalfspacesolution}, which was proved in \cite{AKW} for all $p \in (p_c,p_s)$ and $N \ge 3$ without excluding the locally finite set, by using a different method from ours. The gluing constructions both in \cite{DMP} and \cite{AKW} require a right inverse whose bound is uniform as the concentration scale tends to zero. The locally finite set in \cite{AKW} arises from a non-degeneracy condition to guarantee the existence of such a right inverse. In contrast, our gluing argument allows degenerate connecting solutions, without assuming global uniqueness and non-degeneracy modulo scaling. Instead, we use the nonzero degree of the complete compact set and the explicit error estimates for an exact inner-outer system. We emphasize that this allows us to cover all subcritical exponents in $(p_c,p_s)$ for $N \ge 2$, while the lack of non-degeneracy prevents us from addressing the case when the singular set is a general boundary submanifold.
   
   \section{The half-space case: proof of Proposition \ref{prophalfspacesolution}}
   
   We look for solutions of \eqref{eqhalfspace} of the form
   \[
   u(x) = |x|^{-\frac{2}{p-1}}\phi(-\log |x|, z),
   \]
   where $z = \frac{x}{|x|}$, so that the equation $-\Delta u = u^p$ reads in terms of the function $\phi$ defined for $t \in \R$ and $z \in \mathbb{S}^{N-1}_+$, as
   \begin{align} \label{eq:phi(t,x)}
   	\phi_{tt} + a_p \phi_t + \Delta_{\mathbb{S}^{N-1}}\phi - b_p \phi + \phi^p = 0 \quad \text{in } \C := \R \times \mathbb{S}^{N-1}_+; \qquad \phi = 0 \quad \text{on } \R \times \partial\mathbb{S}^{N-1}_+.
   \end{align}
   The principal construction is a bounded positive solution satisfying
   \begin{align} \label{con:connection}
   	\phi(t,\cdot) \to 0 ~~ \text{ as } t \to -\infty, \quad \phi(t,\cdot) \to \phi_p ~~ \text{ as } t \to +\infty;
   \end{align}
   such a solution is called a {\it connection}.
   
   \subsection{The angular linearization} \label{subangular}
   
   We consider the linearization of problem \eqref{eqhalfsphere} at $\phi = \phi_p$ and all spectral statements in this subsection concern the real Hilbert space $L^2_{\text{ax}}(\s)$ composed of axial $L^2$ functions. For this, we introduce the linear operator:
   \begin{align*}
   	\A_p = -\Delta_{\mathbb{S}^{N-1}} + b_p - p\phi_p^{p-1}, \quad \text{Dom}(\A_p) = H^2_{\text{ax}}(\s) \cap H^1_{0,\text{ax}}(\s). 
   \end{align*}
   The operator $\A_p$ is self-adjoint with compact resolvent, and its quadratic form, on $H^1_{0,\text{ax}}(\s)$, is 
   \begin{align*}
   	Q_p[w] := \int_{\s}\left(|\nabla_{\mathbb{S}^{N-1}}w|^2 + b_pw^2 - p\phi_p^{p-1}w^2\right)d\sigma.
   \end{align*}
   Let $\la_1(p)$ and $\la_2(p)$ denote the first and second axial eigenvalues of $\A_p$. We have the following result.
   \begin{proposition}\label{propaxialnondegenerate}
   	Let $p \in (p_c,p_s)$. The axial eigenvalues of $\A_p$ are simple and satisfy
   	\begin{align*}
   		\la_1(p) < 0 < \la_2(p).
   	\end{align*}
    On any compact $K \Subset (p_c,p_s)$, $p \mapsto \phi_p$ is $C^1$ into $W^{2,q}(\s)$ for every $q > 1$, and we have
    \begin{align} \label{eqinfpinK>0}
    	\inf_{p \in K}\min\{-\la_1(p), \la_2(p)\} > 0.
    \end{align}
   \end{proposition}

   \begin{proof}
   	Repeating arguments in the proof of \cite[Lemma 2.1]{WX}, we know the axial eigenvalues of $\A_p$ are simple and $\la_1(p) < 0 \le \la_2(p)$. The axial non-degeneracy has been already obtained in Shioji-Watanabe \cite[Theorem 12]{SW}. Hence, $\la_2(p) > 0$. Moreover, the map
   	\[
   	(p,w) \mapsto -\Delta_{\mathbb{S}^{N-1}}w + b_p w - (w_+)^{p}
   	\]
   	is $C^1$ from a neighborhood in  $\R\times X_q$ to $L^q_{\text{ax}}$ when $q>N-1$, where $X_q = W^{2,q}_{\text{ax}} \cap W^{1,q}_{0,\text{ax}}$. The derivative in $w$ at $\phi_p$ is an isomorphism by the trivial kernel and elliptic Fredholm theory on the compact hemisphere. Then applying implicit function theorem we obtain the $C^1$ dependence. When $1 < q \le N-1$, take a finite exponent $r > \max\{q,N-1\}$. By the inclusion $X_r \hookrightarrow X_q$ we can get the $C^1$ dependence. Finally, the inequality \eqref{eqinfpinK>0} follows from the eigenvalue continuity and compactness of $K$. The proof is complete.
   \end{proof}
   	
   	Hopf's lemma together with Proposition~\ref{propaxialnondegenerate} implies, for $p$ in any compact interval $K\Subset(p_c,p_s)$,
   	\begin{equation}\label{eq:phidistance}
   		c_K\text{dist}(z,\partial \s)\le\phi_p(z)\le C_K\text{dist}(z,\partial \s),
   		\qquad |\partial_p\phi_p(z)|\le C_K\text{dist}(z,\partial \s).
   	\end{equation}
   	These bounds will also control differentiation with respect to the exponent at the lateral boundary.
   	
   \subsection{A priori bounds and limiting states} \label{subprioribounds}
   
   We start this subsection with a priori bounds for bounded nonnegative classical solutions of \eqref{eq:phi(t,x)}.
   
   \begin{lemma}\label{lem:bound}
   	Let $J\Subset(1,p_s)$. There is a constant $M_J$ such that every bounded nonnegative classical solution of \eqref{eq:phi(t,x)}, with $p\in J$, satisfies
   	\[
   	\|\phi\|_{L^\infty(\C)}\le M_J.
   	\]
   	Its $C^{2,\beta}$ norms on all unit cylindrical slabs, up to the lateral boundary, are uniformly bounded for a sufficiently small fixed $\beta>0$.
   \end{lemma}
   
   \begin{proof}
   	Arguing by contradiction, we suppose the existence of sequences $\{p_j\} \subset J$ and $\{\phi_j\}$ such that $M_j=\|\phi_j\|_{L^\infty(\C)}\to\infty$, with $p_j\to p_*\in J$. Choose $(t_j,z_j) \in \C$ such that $\phi_j(t_j,z_j) \ge M_j/2$. Let
   	\[
   	\varepsilon_j=M_j^{-(p_j-1)/2} \to 0.
   	\]
   	Choose an orthonormal basis $e_{j,1},\cdots,e_{j,N-1}$ in the tangent space
   	\[
   	T_{z_j}\mathbb{S}^{N-1} = \{\zeta \in \R^N: \zeta \cdot z_j = 0\}
   	\]
   	with
	\[
	e_{j,N-1} = \frac{e_N - (z_j)_N z_j}{\sqrt{1-(z_j)_N^2}} \quad (z_j \neq e_N),
	\]
	where $e_N = (0,\cdots,0,1)$ is the north pole, and denote 
   	\[
   	E_j \xi = \sum_{i=1}^{N-1}\xi_i e_{j,i}, \quad \xi \in \R^{N-1}.
   	\]
   	Now we define
   	\[
   	\Phi_j(\xi) = \frac{z_j+E_j\xi}{\sqrt{1+|\xi|^2}}.
   	\]
   	Since $z_j \perp E_j\xi$,
   	\[
   	|z_j+E_j\xi|^2 = 1 + |\xi|^2,
   	\]
   	so $\Phi_j(\xi) \in \mathbb{S}^{N-1}$ and $\Phi_j(0) = z_j$. Set
   	\[
   	w_j(s,y) = \frac{1}{M_j} \phi_j(t_j + \varepsilon_j s,\Phi_j(\varepsilon_j y)), \quad s \in \R, ~~ y \in \R^{N-1} ~\text{satisfying}~ \Phi_j(\varepsilon_j y) \in \s.
   	\]
   	It is clear that $0 \le w_j \le 1$ and $w_j(0,0) = \phi_j(t_j,z_j)/M_j \ge 1/2$, and direct computations yield
   	\[
   	(w_j)_{ss} + \Delta_{g_j} w_j + \varepsilon_j a_{p_j}(w_j)_s - \varepsilon_j^2 b_{p_j} w_j + w_j^{p_j} = 0, \quad s \in \R, ~~ y \in \R^{N-1}, ~ y_{N-1} > -\delta_j := - \frac{(z_j)_N}{\varepsilon_j\sqrt{1-|(z_j)_N|^2}},
   	\]
   	and $w_j = 0$ on $\{y_{N-1} = -\delta_j\}$, where $\Delta_{g_j}$ denotes the Laplace-Beltrami operator with respect to the variable $y$ given as
	\[
	\left(1 + \varepsilon_j^2|y|^2\right)\left[(\delta_{ab}+\varepsilon_j^2y_ay_b)\partial_{y_ay_b} + 2\varepsilon_j^2y \cdot \nabla_y\right].
	\]
	The coefficients of the rescaled angular operator converge locally uniformly to those of $\Delta_y$. When $z_j$ is exactly the north pole, $\delta_j = +\infty$. If $\delta_j\to0$, uniform local boundary gradient estimates yield $w_j(0,0)\le C\delta_j\to0$, contradicting $w_j(0,0)\ge1/2$. Letting $j \to \infty$ we obtain a nonnegative limit solving
   	\[
   	-\Delta w=w^{p_*}
   	\]
   	either in $\R^N$ or in a half-space with zero boundary condition. The whole-space Liouville theorem \cite[Theorem~6.1]{GS} and the half-space theorem \cite[Section~3, Theorem~2 and Remark~2]{Dancer} exclude both limits. For $N=2$ the whole-space conclusion also follows directly from the nonexistence of a positive nonconstant superharmonic function bounded below on $\R^2$. This proves the uniform bound. Standard local and boundary elliptic estimates, first in $W^{2,q}$ and then in H\"older spaces, give the last assertion.
   \end{proof}

   For a bounded nonnegative solution, we define
   \begin{equation}\label{eq:energy}
   	\cE_p[\phi](t)=\int_{\s}\left(
   	\frac12\phi_t^2-\frac12|\nabla_{\mathbb{S}^{N-1}}\phi|^2-\frac{b_p}{2}\phi^2
   	+\frac1{p+1}\phi^{p+1}\right)d\sigma.
   \end{equation}
   Multiplying \eqref{eq:phi(t,x)} by $\phi_t$ yields
   \begin{equation}\label{eq:energydiss}
   	\frac{\mathrm d}{\mathrm dt}\cE_p[\phi](t)
   	=-a_p\int_{\s}\phi_t^2d\sigma,
   \end{equation}
   where we have used $\phi_t=0$ on the lateral boundary.
   
   \begin{lemma}\label{lem:classification}
   	A bounded nonnegative solution of \eqref{eq:phi(t,x)} is either $0$, $\phi_p$, or a positive connection. The limits of a connection hold in $C^2(\overline \s)$, its energy satisfies
   	\begin{equation}\label{eq:energyrange}
   		-J_p\le\cE_p[\phi](t)\le 0,
   		\qquad
   		J_p=\frac{p-1}{2(p+1)}\int_{\s}\phi_p^{p+1}>0,
   	\end{equation}
   	and
   	\begin{equation}\label{eq:dissipationtotal}
   		a_p\int_\R\int_{\s}\phi_t^2d\sigma d t=J_p.
   	\end{equation}
   \end{lemma}
   
   \begin{proof}
   	Let $\phi$ be a bounded nonnegative solution of \eqref{eq:phi(t,x)}. By Lemma~\ref{lem:bound} and \eqref{eq:energydiss}, the energy $\cE_p[\phi](t)$ is bounded, and nonincreasing in $t \in \R$. Hence $\phi_t\in L^2(\C)$. For a sequence $t_j\to+\infty$ or $t_j\to-\infty$, translate $\phi$ by $t_j$ and pass to a subsequence in local $C^2$. The integral of $\phi_t^2$ over any fixed translated slab tends to zero. The limit is therefore independent of $t$ and is either $0$ or the unique solution $\phi_p$ of \eqref{eqhalfsphere}.
   	
   	Note that
   	\[
   	\cE_p[0]=0,\qquad \cE_p[\phi_p]=-J_p.
   	\]
   	At either end the energy already has a single limit, so the two different equilibria cannot both occur as subsequential limits at that end. Thus each end has a unique limit. If the limits coincide, \eqref{eq:energydiss} forces $\phi_t\equiv0$. Otherwise the direction is necessarily $0$ at $-\infty$ and $\phi_p$ at $+\infty$. Positivity of the connection follows from the strong maximum principle. Finally, \eqref{eq:dissipationtotal} follows by the integration of \eqref{eq:energydiss}. The proof is complete.
   \end{proof}
 
   Define the linear phase functional
   \begin{equation}\label{eq:phase}
   	\ell(\phi)=\int_{\s}\left(\phi(0,z)+\frac1N \phi_t(0,z)\right)\varphi_1(z)d\sigma,
   	\quad
   	m_p=\int_{\s}\phi_p\varphi_1 d\sigma,
   \end{equation}
   where $\varphi_1$ is given in \eqref{eq:phi1}. A connection is \textit{normalized} if
   \begin{equation}\label{eq:normalization}
   	\ell(\phi)=\tfrac12m_p.
   \end{equation}
   Every connection has a translate satisfying this condition: the same expression at time $s$ tends to zero as $s\to-\infty$ and to $m_p$ as $s\to+\infty$. We do not assume that the normalizing translate is globally unique.
   
   \begin{lemma}\label{lem:uniformends}
   	Let $\{p_j\} \subset K\Subset(p_c,p_s)$ satisfy $p_j\to p_*$ and $\{\phi_j\}$ be a sequence of normalized axial connections. Then, up to a subsequence there is a limit $\phi_*$ such that $\phi_j \to \phi_*$ in $C^2$ on local cylindrical slabs up to the lateral boundary, and $\phi_*$ is a normalized axial connection at $p = p_*$. Moreover, the convergence of $\phi_j$ to $0$ at $-\infty$ and to $\phi_{p_j}$ at $+\infty$ is uniform with respect to $j$ in the local $C^2$ slab norms.
   \end{lemma}
   
   \begin{proof}
   	Lemma~\ref{lem:bound} yields the existence of the limit $\phi_*$. Moreover, using the phase condition $\ell(\phi_j) = 1/2m_{p_j}$ one has $\ell(\phi_*) = 1/2m_{p_*}$. This shows that $\phi_*$ is neither $0$ nor $\phi_{p_*}$. Then, by Lemma~\ref{lem:classification} we get that $\phi_*$ is a normalized connection.
   	
   	To obtain the uniform convergence as $t \to +\infty$, suppose by contradiction the existence of a subsequence of $\{t_j\}$, also denoted by $t_j$, such that  $t_j\to+\infty$ and that $\phi_j(t_j,\cdot)$ does not approach $\phi_{p_j}$. Note that $\phi_j(t+t_j,z)$ converges locally to a bounded nonnegative solution $V$ up to a subsequence. For any fixed $T$ and bounded $s$, for sufficiently large $j$,
   	\[
   	-J_{p_j}\le\cE_{p_j}[\phi_j](t_j+s)
   	\le\cE_{p_j}[\phi_j](T).
   	\]
   	First let $j\to\infty$ and then $T\to+\infty$. Since $\phi_*$ is a connection, the right-hand side tends to $-J_{p_*}$. Hence, $\cE_{p_*}[V]\equiv-J_{p_*}$, so \eqref{eq:energydiss} implies $V\equiv\phi_{p_*}$, a contradiction. For the uniform convergence as $t_j \to -\infty$, repeating above arguments by replacing $\phi_{p_j}$ by $0$ and using for fixed $T>0$, bounded $s$ and large $j$
   	\[
   	0\ge\cE_{p_j}[\phi_j](t_j+s)
   	\ge\cE_{p_j}[\phi_j](-T),
   	\]
   	we can get a desired contradiction to prove the result. Finally, by the elliptic compactness we can improve the uniform convergence in the local $C^2$ slab norms and complete the proof.
   \end{proof}
   
   \subsection{Compactness in a fixed weighted space} \label{subcompactness}

   \begin{lemma}\label{lem:green}
   	Suppose $a>0$. Let $V$ belong to $L^\infty_{\text{ax}}(\s)$ and let $B=-\Delta_{\mathbb{S}^{N-1}}+V(z)$ act on $L^2_{\text{ax}}(\s)$ with form domain $H^1_{0,\text{ax}}(\s)$. Assume that $0\notin\sigma(B)$. On bounded $L^2_{\text{ax}}(\s)$-valued functions, the operator
   	\[
   	L=\partial_{tt}+a\partial_t-B
   	\]
   	has a unique bounded inverse given by convolution with an operator kernel $G_B$. Moreover, there exist $C,\kappa>0$ such that
   	\begin{equation}\label{eq:greenbound}
   		\|G_B(t)\|_{L^2\to L^2}
   		+\|G_B'(t)\|_{L^2\to L^2}
   		\le C e^{-\kappa|t|}\qquad(t\ne0).
   	\end{equation}
   	Constants $C$ and $\kappa$ can be chosen uniformly in a compact family for which $a$ stays positive and the spectra of $B$ stay away from zero.
   \end{lemma}
   
   \begin{proof}
   	The operator $B$ is self-adjoint, bounded below, and has compact resolvent. There exists an orthonormal basis $\{e_j\}$ of $L^2_{\text{ax}}(\s)$ such that
   	\[
   	B e_j = \lambda_j e_j, \quad \lambda_1 \le \lambda_2 \le \cdots, \quad \lambda_j \to +\infty,
   	\]
   	where $\{\lambda_j\}$ is the sequence of eigenvalues of $B$. The number of negative eigenvalues is finite and we assume that $\la_{m-1} < 0 < \la_m$, where $m = 1$ means that there is no negative eigenvalue. For $L^2_{\text{ax}}(\s)$-valued functions $u, f$, we can write
   	\[
   	u(t) = \sum_{j=1}^\infty u_j(t)e_j, \quad f(t) = \sum_{j=1}^\infty f_j(t)e_j,
   	\]
   	where
   	\[
   	u_j(t) = \langle u(t),e_j \rangle_{L^2}, \quad f_j(t) = \langle f(t),e_j \rangle_{L^2}.
   	\]
   	Projecting $Lu = f$ onto $\text{span}\{e_j\}$, we have
   	\[
   	u''_j + a u'_j - \lambda_j u_j = f_j.
   	\]
   	
   	For an eigenvalue $\lambda$ of $B$, let $r_\pm$ be the roots of
   	\begin{equation}\label{eq:roots}
   		r^2+ar-\lambda=0,
   	\end{equation}
    given by
    \[
    r_\pm = \frac{-a \pm \sqrt{a^2 + 4\la}}{2}.
    \]
   	If $\lambda>0$, then $\la \ge \la_m$ and $r_-<0<r_+$. We define
   	\[
   	g_\lambda(t)=-\frac1{r_+-r_-}
   	\begin{cases}e^{r_-t},&t\ge0,\\ e^{r_+t},&t<0.\end{cases}
   	\]
   	Let $\gamma_\lambda = r_+ > 0$. Note that $-r_- = a + r_+ > \gamma_\la$. Since $\gamma_\la$ is increasing with respect to $\la$, we have
    \[
    \gamma_\la \ge \gamma_+ := \gamma_{\la_m}.
    \]
    Computing directly, we obtain
    \begin{align} \label{eagamma+}
    	|g_\la(t)| \le \frac{1}{\sqrt{a^2+4\la}}e^{-\gamma_\la|t|} \le \frac{1}{\sqrt{a^2+4\la_m}}e^{-\gamma_+|t|}, \quad |g'_\la(t)| \le e^{-\gamma_\la|t|} \le e^{-\gamma_+|t|}.
    \end{align}

   	If $\lambda<0$, both roots have negative real parts and we define for $r_- \neq r_+$
   	\[
   	g_\lambda(t)=\begin{cases}
   		(e^{r_+t}-e^{r_-t})/(r_+-r_-),&t\ge0,\\0,&t<0.
   	\end{cases}
   	\]
   	If $r_-=r_+$, $g_\lambda(t) = 0$ for $t < 0$ and $g_\lambda(t) = te^{-at/2}$ for $t \ge 0$. A direct computation yields that $g_\lambda(t)$ is always real for complex conjugate roots, which occur when $a^2 + 4\lambda < 0$. Setting
   	\begin{equation}
   	\gamma_\la = 
   	\left\{
   	\begin{aligned}
   		& \frac{a}{2} - \sqrt{\frac{a^2}{4} + \la} , & \quad \text{if } a^2 + 4\lambda > 0, \\
   		& \frac{a}{2}, & \quad \text{if } a^2 + 4\lambda \le 0,
   	\end{aligned}   	
   	\right.
   	\end{equation}
   	it follows direct computations that
   	\begin{align*}
   		|g_\la(t)| \le t e^{-\gamma_\la t}, \quad |g'_\la(t)| \le \left(1 + \frac{a}{2}t\right) e^{-\gamma_\la t}, \quad t > 0.
   	\end{align*}
   Let $\gamma_- := \min_{\la_j < 0}\gamma_{\la_j} > 0$ and we choose
   \[
   0 < \kappa < \gamma_-. 
   \]
   We further have 
   \begin{align} \label{eagamma-}
   	|g_\la(t)| \le \frac{1}{e(\gamma_\la - \kappa)} e^{-\kappa t}, \quad |g'_\la(t)| \le \left(1 + \frac{a}{2e(\gamma_\la - \kappa)}\right) e^{-\kappa t}.
   \end{align}

  Decrease $\kappa$ if necessary so that it also holds $\kappa < \gamma_+$. Combining \eqref{eagamma+} and \eqref{eagamma-} yields
  \begin{align*}
  	|g_{\la_j}(t)| \le C_0 e^{-\kappa |t|}, \quad |g'_{\la_j}(t)| \le C_1 e^{-\kappa |t|},
  \end{align*}
  for two positive constants $C_0, C_1$ independent of $j$. Define the operator kernel by
  \begin{align*}
  	G_B(t)h = \sum_{j=1}^\infty g_{\la_j}(t)\langle h,e_j \rangle_{L^2} e_j.
  \end{align*}
  Then \eqref{eq:greenbound} follows with $C= C_0 + C_1$. It is readily seen that there is no nonzero bounded solution $u$ solving $Lu = 0$ and the convolution with $G_B$ is the unique bounded inverse. Finally, checking the proof we can obtain the uniformity of $C$ and $\kappa$ once $a$ stays positive and the spectra of $B$ stay away from zero. The proof is complete.
   \end{proof}

   \begin{lemma}\label{lem:exptails}
   	For every compact $K\Subset(p_c,p_s)$, there exist $T,C,\eta>0$ such that every normalized axial connection with $p\in K$ satisfies
   	\begin{align}
   		&\|\phi\|_{W^{2,q}((t-1,t+1)\times \s)}
   		\le Ce^{\eta t} &&(t\le-T),\label{eq:tailminus}\\
   		&\|\phi-\phi_p\|_{W^{2,q}((t-1,t+1)\times \s)}
   		\le Ce^{-\eta t} &&(t\ge T),\label{eq:tailplus}
   	\end{align}
   	for $q>1$, where the constant $C$ also depends on $q$. The estimates also hold in $C^1$ on these slabs. For a single connection at its negative end, \eqref{eq:tailminus} holds with any
   	\begin{equation}\label{eq:allowedrate}
   		0<\eta<N-1-\frac{2}{p-1}.
   	\end{equation}
   \end{lemma}
   
   \begin{proof}
   	Set
   	\begin{equation}\label{eq:endoperators}
   		B_{-,p}=-\Delta_{\mathbb{S}^{N-1}}+b_p>0,\qquad B_{+,p}=\A_p.
   	\end{equation}
   	By Proposition~\ref{propaxialnondegenerate}, $B_{+,p}$ has one negative eigenvalue and no zero eigenvalue. Lemma~\ref{lem:green} applies at both ends with a uniform rate dependent on $K$.
   	
   	Consider the positive end and write $h=\phi-\phi_p$, which satisfies
   	\begin{align} \label{eqwitherrorrho}
		(\partial_{tt}+a_p\partial_t-B_{+,p})h=\rho(t,z)h,
	\end{align}
   	where $\rho\to0$ as $t\to+\infty$, uniformly over the normalized family. This follows from Lemma~\ref{lem:uniformends} and the uniform continuity of the derivative of $s\mapsto(s_+)^p$ on bounded intervals, uniformly for $p\in K$. Take a smooth cut-off function $\chi(t)$ such that $\chi(t) = 0$ for $t \le T$ and $\chi(t) = 1$ for $t \ge T+1$ and let $\widetilde h = \chi h$. Then, on the entire cylinder
   	\begin{equation}\label{eq:cutoffeq}
   		(\partial_{tt}+a_p\partial_t-B_{+,p})\widetilde h
   		=\rho\widetilde h+g, \quad g = 2\chi' h_t + (\chi'' + a_p \chi')h.
   	\end{equation}
	Note that $\sup_t \|g(t)\|_{L^2}$ is uniformly bounded on the compact $K$.
   	
   	To avoid presupposing weighted boundedness, we use
   	\[
   	w_R(t)=\exp\big(\eta\min\{(t-T)_+,R\}\big),\qquad R>0.
   	\]
   	They satisfy $w_R(t)/w_R(s)\le e^{\eta|t-s|}$. Convolve \eqref{eq:cutoffeq} with the Green kernel and take the weighted supremum in $L^2(\s)$. If $0<\eta<\kappa$, \eqref{eq:greenbound} yields
   	\[
   	\sup_t w_R(t)\|\widetilde h(t)\|_{L^2}
   	\le C_\eta \|\rho\|_{L^\infty([T,\infty)\times \s)}\sup_t w_R(t)\|\widetilde h(t)\|_{L^2}
   	+C_\eta\sup_t w_R(t)\|g(t)\|_{L^2},
   	\]
	that is
	\[
	 \left(1-C_\eta \|\rho\|_{L^\infty([T,\infty)\times \s)}\right) \sup_t w_R(t)\|\widetilde h(t)\|_{L^2} \le C_\eta\sup_t w_R(t)\|g(t)\|_{L^2}.
	\]
   	Choose $T$ large so that the coefficient of the left side is positive. Moreover, since $g(t) \equiv 0$ for $t$ outside $[T,T+1]$ and on $[T,T+1]$, we have $w_R(t) \le e^\eta$. Hence, 
	\[
	\sup_t w_R(t)\|g(t)\|_{L^2} = \sup_{T \le t \le T+1} w_R(t)\|g(t)\|_{L^2} \le e^\eta\sup_t\|g(t)\|_{L^2}.
	\]
	Letting $R\to\infty$, for $t \ge T+1$ we deduce that $e^{\eta(t-T)}\|h(t)\|_{L^2}$ is uniformly bounded. Using the estimate for $G_B'$, we can further derive the uniform bound for $e^{\eta(t-T)}\|h_t(t)\|_{L^2}$ for $t \ge T + 1$. In view of \eqref{eqwitherrorrho}, standard bootstrap arguments yield \eqref{eq:tailplus} and the $C^1$ estimate. Similar arguments give a proof for negative end.
   	
   	For the last assertion, the smallest positive root of \eqref{eq:roots} for $B_{-,p}$ is $N-1- 2/(p-1)$, whereas the absolute value of its corresponding negative root is $2/(p-1)+1$. The former is smaller for $p \in K$. Thus the negative-end Green estimate can be used with every rate in \eqref{eq:allowedrate}.
   \end{proof}
   
   Fix $q>N$. Let $\rho\in C^\infty(\R)$ be nonnegative and equal to $|t|$ for $|t|\ge2$. For $\delta>0$ define the real Banach spaces
   \begin{align}
   	\cX_\delta={}&\{w\in W^{2,q}_{\mathrm{loc}}(\overline\C):
   	w\text{ is axial},\ {\rm Tr}_{\R\times\partial \s}w=0,
   	\ e^{\delta\rho}w\in W^{2,q}(\C)\},\label{eq:Xspace}\\
   	\cY_\delta={}&\{f\in L^q_{\mathrm{loc}}(\C):f\text{ is axial},
   	\ e^{\delta\rho}f\in L^q(\C)\}.
   	\label{eq:Yspace}
   \end{align}
   Their norms are $\|e^{\delta\rho}w\|_{W^{2,q}}$ and $\|e^{\delta\rho}f\|_{L^q}$, respectively. In particular, elements of $\cX_\delta$ and their first derivatives decay at both ends. Take a cut-off function $\chi\in C^\infty(\R)$ such that $0\le\chi\le1$, $\chi=0$ on $(-\infty,-1]$ and $\chi=1$ on $[1,\infty)$.
   
   \begin{proposition}\label{prop:compact}
   	For every compact $K\Subset(p_c,p_s)$ there is $\delta_K>0$ such that, for $0<\delta<\delta_K$, the set
   	\begin{equation}\label{eq:compactset}
   		\{(p,\phi-\chi\phi_p):p\in K,\ \phi\text{ is a normalized axial connection}\}
   	\end{equation}
   	is compact in $K\times\cX_\delta$.
   \end{proposition}
   
   \begin{proof}
   	Choose $\eta$ from Lemma~\ref{lem:exptails} and take $0<\delta<\eta$. Given a sequence $\{p_j\} \subset K$ and a connection sequence $\{\phi_j\}$, up to a subsequence, we have $p_j \to p_* \in K$ and Lemma~\ref{lem:uniformends} yields the existence of $\phi_*$ such that $\phi_j \to \phi_*$ in $C^2_{{\rm loc}}(\R \times \overline{\s})$ up to the lateral boundary. Moreover, by Proposition~\ref{propaxialnondegenerate}, $\phi_{p_j} \to \phi_{p_*}$ in $W^{2,q}(\s)$. Set
	\[
	u_j = \phi_j - \phi_*, \quad P_j = \partial_{tt} + a_{p_j}\partial_t + \Delta_{\cs} - b_{p_j}.
	\]
	Note that
	\[
	P_j u_j = -(a_{p_j}-a_{p_*})(\phi_*)_t + (b_{p_j} - b_{p_*})\phi_* - \phi_j^{p_j} + \phi_*^{p_*}, \quad u_j = 0 ~\text{ on } \R \times \partial \s.
 	\]
	On any fixed finite cylinder $(-R,R) \times \s$, by standard elliptic estimates we obtain $u_j \to 0$ in $W^{2,q}$. Let $w_j = \phi_j - \chi\phi_{p_j}$ and $w_* = \phi_* - \chi\phi_{p_*}$. On any fixed finite cylinder, $w_j - w_*$ converges to $0$ in $W^{2,q}$. Now we define
	\[
	Q_m^+ = (m,m+1) \times \s, \quad Q_m^- = (-m-1,-m) \times \s,
	\]
	where $m > 0$ is a sufficiently large integer. Using Lemma~\ref{lem:exptails} we deduce that 
	\[
	\|w_j\|_{W^{2,q}(Q_m^\pm)} \le Ce^{-\eta m}, \quad \|w_*\|_{W^{2,q}(Q_m^\pm)} \le Ce^{-\eta m}.
	\]
	By the definition of $\rho$ and a direct computation, we get
	\[
	\|e^{\delta\rho}w_j\|_{W^{2,q}(Q_m^\pm)} \le Ce^{-(\eta-\delta)m}, \quad \|e^{\delta\rho}w_*\|_{W^{2,q}(Q_m^\pm)} \le Ce^{-(\eta-\delta)m}
	\]
	where $C$ may vary. The sum of the $q$th powers of the slab norms from $m = M$ is dominated by
	\[
	\sum_{m=M}^\infty e^{-q(\eta-\delta)m} = \frac{e^{-q(\eta-\delta)M}}{1-e^{-q(\eta-\delta)}}.
	\]
	Letting $M$ large enough first and then using the uniform convergence on $(-M,M)\times \s$ proves convergence in $\cX_\delta$. The proof is complete.
   \end{proof}

   \subsection{A cylindrical Fredholm index calculation} \label{subindex}
   
   \begin{lemma}\label{lem:stepindex}
   	Let $a>0$. Let $B_- = -\Delta_{\mathbb{S}^{N-1}} + V_-$ and $B_+ = -\Delta_{\mathbb{S}^{N-1}} + V_+$ be axial self-adjoint Dirichlet realizations of $-\Delta_{\mathbb{S}^{N-1}}$ plus bounded real potentials. Assume
   	\[
   	B_->0,\qquad
   	\lambda_1(B_+)<0<\lambda_2(B_+),
   	\]
   	and $\lambda_1(B_+)$ is simple. Define
   	\begin{equation}\label{eq:stepoperator}
   		\widehat Lh=h_{tt}+ah_t-B(t)h,
   		\qquad
   		B(t)=\begin{cases}B_-,&t<0,\\B_+,&t>0.\end{cases}
   	\end{equation}
   	For all sufficiently small $\delta>0$, the operator
   	\[
   	\widehat L:\cX_\delta\rightarrow\cY_\delta
   	\]
   	is surjective and has a one-dimensional kernel. In particular, its Fredholm index is one.
   \end{lemma}
   
   \begin{proof}
   	\textit{Step 1: the matching space.}
   	Let $\psi$ be the normalized eigenfunction of $B_+$ corresponding to $\lambda_1(B_+)$ and let $Q$ be the orthogonal projection onto $\text{span}\{\psi\}^{\perp}$. Define
   	\begin{equation}\label{eq:Dminusplus}
   		D_- =\frac{-a+(a^2+4B_-)^{1/2}}2,
   		\quad
   		D_+ =\frac{a+(a^2+4B_+|_{{\rm Ran } Q})^{1/2}}2,
   	\end{equation}
   	where $D_+$ is defined only on the positive spectral subspace ${\rm Ran } Q$. Note that both $D_-$ and $D_+$ are positive operators. The operator $D_-$ has the following form domain:
   	\[
   	E=\text{Dom}(B_-^{1/4})\subset L^2_{\text{ax}}(\s),
   	\qquad \|h\|_E^2=\|B_-^{1/4}h\|_{L^2(\s)}^2.
   	\]
   	Equivalently, $E$ is the interpolation space between $L^2_{\text{ax}}(\s)$ and $H^1_{0,\text{ax}}(\s)$ with interpolation exponent $1/2$. Denote $E^*$ the dual space of $E$. Since $B_-$ and $B_+$ have same second-order part $-\Delta_{\mathbb{S}^{N-1}}$ and the potentials are bounded, $E$ is also the form domain of $D_+Q$. Now we can define the operator
   	\begin{equation}\label{eq:matchingK}
   		\mathcal K=D_-+D_+Q:E\rightarrow E^*
   	\end{equation}
   	by its quadratic form as
	\[
	\langle\mathcal Kh,\cdot\rangle
   		=\langle D_-h,\cdot \rangle + \langle D_+Q h,\cdot \rangle
	\]
	The form is bounded on $E$, and
   	\begin{equation}\label{eq:Kcoercive}
   		\langle\mathcal Kh,h\rangle
   		=\langle D_-h,h \rangle + \langle D_+Q h,Qh \rangle 
   		\ge \ep\|h\|_E^2
   	\end{equation}
   	for some $\ep > 0$, by the spectral assumptions of $B_-, B_+$ and the comparison between the spectrum of $B_-, B_+$ and $D_-, D_+$. Applying Lax-Milgram, we deduce that $\mathcal K$ is an isomorphism from $E$ onto $E^*$.
   	
   	\textit{Step 2: homogeneous solutions.}
   	In this step we would like to solve $\widehat Lh=0$. Since we search for solutions in $\cX_\delta$, a solution $h$ decays as $|t| \to \infty$. Assume that $h$ has value $h_0$ at zero, that is, $h(0)=h_0$. Then, on $t < 0$, $h$ is given by
   	\[
   	h_-(t)=e^{D_-t}h_0,
   	\quad h_-'(0)=D_-h_0,
   	\]
   	and on $t>0$, the positive spectral component of $h$ with the same trace $h_0$ is
   	\[
   	Qh_+(t)=e^{-D_+t}Qh_0,
   	\quad Qh_+'(0)=-D_+Qh_0.
   	\]
   	It remains to determine $h_+$ in the $\psi$ direction. Set $y(t) = \langle h_+(t),\psi \rangle$. Projecting $\widehat Lh=0$ onto $\text{span}\{\psi\}$ yields 
	\[
	y'' + ay' - \la_1(B_+) y = 0.
	\]
	Note that $a>0$ and $\la_1(B_+) < 0$. For any initial value $y(0) = \langle h_0,\psi \rangle$ and its initial derivative $y'(0) = c \in\R$, we can solve a decaying solution $y(t)$.
   	
   	To obtain the solution, we need to match the two derivatives $h_-'(0)$ and $h_+'(0)$, implying
   	\begin{equation}\label{eq:matchhom}
   		\mathcal Kh_0-c\psi=0.
   	\end{equation}
   	Since $\mathcal K$ is an isomorphism, for each $c$ there is exactly one trace $h_0=c\mathcal K^{-1}\psi$. If $c=0$, the solution is zero. Thus the space of decaying homogeneous solutions has dimension one.
   	
   	\textit{Step 3: a right inverse in weighted Hilbert spaces.}
   	First, for $f \in L^2(\C)$ compactly supported in $t$, we search for $u$ such that $\widehat Lu=f$. Convolution with the ordinary Dirichlet half-line Green kernels provides a decaying solution $w_-$ on the negative half-cylinder and positive spectral component $Qw_+$ on $t > 0$  with zero value at $t=0$: 
	\begin{align*}
		& \widehat Lw_- = f ~~ (t < 0), \quad w_-(0) = 0, \quad w_-(t) \to 0 ~~ (t \to -\infty), \\
		& \widehat LQw_+ = Qf ~~ (t > 0), \quad Qw_+(0) = 0, \quad Qw_+(t) \to 0 ~~ (t \to +\infty).
	\end{align*}
	Solving 
	\[
	y'' + ay' - \la_1(B_+)y = \langle f,\psi \rangle
	\]
	with $y(0) = 0$ and $y'(0) = 0$, we obtain a decaying solution $y(t)$. Then, $w_+ = Qw_+ + y(t)\psi$ satisfies
	\[
	\widehat Lw_+ = f ~~ (t > 0), \quad w_+(0) = 0, \quad w_+(t) \to 0 ~~ (t \to +\infty).
	\]
   	
   	Let $d_- = w'_-(0), d_+ = w'_+(0)\in E^*$. Note that $d_- \neq d_+$ in general. Given trace $h_0$ and $c \in \R$, we can get $h_-$ and $h_+$ as in Step 2. We aim to find a solution $u = w + h$, where $w=w_-$, $h = h_-$ for $t < 0$, $w=w_+$, $h = h_+$ for $t > 0$, $w(0)= 0$, $h(0) = h_0$. Matching the derivatives $d_- + h'_-(0)$ and $d_+ + h'_+(0)$ yields
   	\begin{equation}\label{eq:matchinhom}
   		\mathcal Kh_0-c\psi=d_+-d_-.
   	\end{equation}
   	Setting $c=0$ and determining $h_0$ by $\mathcal K^{-1}(d_+-d_-)$, we obtain the desired solution $u$ such that $\widehat Lu=f$ holds in the weak sense. By the elliptic regularity, we can get the local $H^2$ regularity of $u$. 
   	
   	Take small $\kappa>0$ dependent on the spectra of $B_-$ and $B_+$. Let $I_k=(k,k+1)$ and 
	\[
	f = \sum_{k \in \mathbb{Z}}f^{(k)}, \quad f^{(k)}(t,z) = 1_{I_k}(t)f(t,z).
	\] 
	Since $f$ has compact support, only finitely many $f^{(k)}$ are nonzero. Due to that $\widehat{L}$ is linear, $u = \sum u^{(k)} = \sum (w^{(k)} + h^{(k)})$, where $u^{(k)}$ is the solution of $\widehat{L}(u^{(k)}) = f^{(k)}$ by a same way as $u$. Using the estimates of the half-line Green kernels, a direct computation yields 
	\begin{align} \label{w^(k)}
		\|w^{(k)}\|_{L^2(I_j \times \s)} \le Ce^{-\kappa|j-k|}\|f^{(k)}\|_{L^2(\C)}.
	\end{align}
	Using
	\[
	\|d_+^{(k)} - d_-^{(k)}\|_{E^*} \le Ce^{-\kappa|k|}\|f^{(k)}\|_{L^2(\C)},
	\]
	together with the homogeneous extensions of $h^{(k)}$ and the boundedness of $\mathcal K^{-1}$, we further get
	\begin{align} \label{h^(k)}
		\|h^{(k)}\|_{L^2(I_j \times \s)} \le Ce^{-\kappa(|j|+|k|)}\|f^{(k)}\|_{L^2(\C)}.
	\end{align}
	Combining \eqref{w^(k)} and \eqref{h^(k)} yields
	\[
	\|u^{(k)}\|_{L^2(I_j \times \s)} \le Ce^{-\kappa(|j-k|)}\|f^{(k)}\|_{L^2(\C)},
	\]
	where $C$ may vary, and so
   	\begin{equation}\label{eq:slabkernelestimate}
   		\|u\|_{L^2(I_j\times \s)} \le \sum_{k}\|u^{(k)}\|_{L^2(I_j \times \s)} 
   		\le C\sum_{k}e^{-\kappa|j-k|}
   		\|f\|_{L^2(I_k\times \s)}.
   	\end{equation}
	This enables us to obtain for sufficiently small $\delta$
	\begin{align}
		\|e^{\delta\rho}u\|_{L^2(\C)} \le C_\delta\|e^{\delta\rho}f\|_{L^2(\C)}.
	\end{align}
	The local and boundary elliptic estimates yield
	\begin{align}
		\|e^{\delta\rho}u\|_{H^2(\C)} \le C_\delta\|e^{\delta\rho}f\|_{L^2(\C)}.
	\end{align}
	Now we can remove the compact support assumption of $f$ by continuous extension to $\{f \in L^2(\C):\|e^{\delta\rho}f\|_{L^2(\C)} < \infty\}$. 
   	
   	\textit{Step 4: passage to $W^{2,q}$ and $L^q$.}
   	For the fixed $q>N$, by local elliptic estimates and a bootstrap argument, we get
   	\begin{equation}\label{eq:localqestimate}
   		\|u\|_{W^{2,q}(I_j\times \s)}
   		\le C_q\left(
   		\| f\|_{L^q((j-2,j+3)\times \s)}
   		+\| u\|_{L^2((j-2,j+3)\times \s)}\right).
   	\end{equation}
   	Then, multiplying \eqref{eq:slabkernelestimate} and \eqref{eq:localqestimate} by $e^{\delta|j|}$ with $\delta > 0$ small enough, and using discrete Young's inequality for $e^{-(\kappa-\delta)|j-k|}$ yields that
   	\begin{equation}\label{eq:rightinversebound}
   		\|u\|_{\cX_\delta}\le C\| f\|_{\cY_\delta}.
   	\end{equation}
   	Therefore, for each $f \in \cY_\delta$, we can find $u \in \cX_\delta$ such that $\widehat{L}u = f$ when $\delta > 0$ is small.
   	
   	Finally, it is straightforward to check that the kernel of $\widehat L:\cX_\delta\rightarrow\cY_\delta$ is contained in the one-dimensional Hilbert-space kernel found in Step 2 and conversely, the homogeneous solution established in Step 2 indeed belongs to $\cX_\delta$. Thus the kernel dimension is exactly one. The proof is complete.
   \end{proof}
   
   \begin{proposition}\label{prop:variableindex}
   	Fix a compact $K\Subset(p_c,p_s)$ and take $\delta>0$ sufficiently small, uniformly for $p\in K$. For every $p\in K$ and $w\in\cX_\delta$, let $\phi=\chi\phi_p+w$. Then, the operator
   	\begin{equation}\label{eq:actualL}
   		L_{p,\phi}=\partial_{tt}+a_p \partial_t+\Delta_{\mathbb{S}^{N-1}}-b_p
   		+p(\phi_+)^{p-1} : \cX_\delta\rightarrow\cY_\delta
   	\end{equation}
   	is Fredholm of index one. Consequently,
   	\begin{equation}\label{eq:phaseindex}
   		h\mapsto (L_{p,\phi}h,\ell(h))
   		:\cX_\delta\rightarrow\cY_\delta\times\R
   	\end{equation}
   	is Fredholm of index zero.
   \end{proposition}
   
   \begin{proof}
   	 We apply Lemma~\ref{lem:stepindex} for the operator $\widehat L$ with $B_- = B_{-,p}$, $B_+ = B_{+,p}$, where $B_{-,p}$ and $B_{+,p}$ are given in \eqref{eq:endoperators}. By the definition given in \eqref{eq:actualL}, the difference between $\widehat L$ and $L_{p,\phi}$ is multiplication by a bounded function that tends uniformly to $0$ as $t \to -\infty$ or $t \to +\infty$. Note that such a multiplication is compact from $\cX_\delta$ to $\cY_\delta$. Hence, $L_{p,\phi}$ has the same Fredholm index as $\widehat L$, namely one. It is clear that $h \mapsto (L_{p,\phi}h,0)$ has Fredholm index zero, since the dimension of kernel is unchanged while the codimension is incremented by one. The map \eqref{eq:phaseindex} is a finite-rank perturbation of $h \mapsto (L_{p,\phi}h,0)$. Thus the index is unchanged. The proof is complete.
   \end{proof}

   \subsection{Proof of Proposition \ref{prophalfspacesolution}} \label{subproof1}
   
   Fix $K=[p_0,P]\Subset(p_c,p_s)$ and a weight $\delta>0$ small enough for both Proposition~\ref{prop:compact} and Proposition~\ref{prop:variableindex}. Define
   \begin{equation}\label{eq:Fmap}
   	\cF: K\times\cX_\delta \to \cY_\delta\times\R; \quad \cF(p,w)=
   	\left(
   	\phi_{tt}+a_p\phi_t+\Delta_{\mathbb{S}^{N-1}}\phi-b_p\phi+(\phi_+)^p,
   	\ \ell(\phi)-\tfrac12m_p
   	\right),\quad \phi=\chi\phi_p+w.
   \end{equation}
   
   \begin{lemma}\label{lem:Fproperties}
   	The map $\cF$ is $C^1$ and $D_w\cF$ is exactly \eqref{eq:phaseindex}. Every zero of $\cF$ corresponds to a normalized axial connection, conversely, every such connection gives a zero of $\cF$. Moreover, the zero set $\{\cF(p,w) = 0\}$ is compact.
   \end{lemma}
   
   \begin{proof}
   	We expand the first component of $\cF$ as
   	\begin{align*}
   		&w_{tt}+a_pw_t+\Delta_{\mathbb{S}^{N-1}}w -b_pw
   		+(\chi''+a_p\chi')\phi_p+(\chi^p-\chi)\phi_p^p\\
   		&\hspace{35mm}+\big[(\chi\phi_p+w)_+^p-(\chi\phi_p)^p\big].
   	\end{align*}  
    Then, it is straightforward to check that it belongs to $\cY_\delta$ when $w \in \cX_\delta$, so the map $\cF$ is well-defined. It is not difficult to see that this map is $C^0$, its derivative with respect to $w$ is \eqref{eq:phaseindex}, and $D_w\cF$ is continuous. It remains to show that
    \[
    D_p\cF = (H_p(w), \ell(\chi\partial_p\phi_p)-\frac12m_p')
    \]
    is continuous, where
    \[
    H_p(w) = a_p'w_t - b_p'w + a_p'\chi'\phi_p + (\chi''+a_p\chi')\partial_p \phi_p + [(\phi_+)^p\log(\phi_+)-\chi\phi_p^p\log\phi_p] + p\chi[(\phi_+)^{p-1}-\phi_p^{p-1}]\partial_p\phi_p,
    \] 
    $(\phi_+)^p\log(\phi_+) = 0$ if $\phi \le 0$, and $a_p', b_p', m_p'$ denote the derivatives with respect to $p$. By Proposition~\ref{propaxialnondegenerate}, $\partial_p\phi_p$ is well-defined. On a finite cylinder, the scalar map $(p,s)\mapsto(s_+)^p$ is $C^1$ since $p>1$, with $\partial_p(s_+)^p=(s_+)^p\log(s_+)$. On the negative end, $\chi = 0$ and $\phi = w$. On the positive end, $\chi = 1$, $\phi = \phi_p + w$, and we focus on terms
    \begin{align} \label{twoterms}
    	[(\phi_+)^p\log(\phi_+)-\phi_p^p\log\phi_p] + p[(\phi_+)^{p-1}-\phi_p^{p-1}]\partial_p\phi_p.
    \end{align}
    By \eqref{eq:phidistance} and the $C^1$ weighted embedding, we have
   	\[
   	|w(t,z)|\le C e^{-\delta t}\text{dist}(z,\partial \s),\quad
   	|\partial_p\phi_p(z)|\le C\phi_p(z).
   	\]
   	Thus, in a neighborhood of any fixed $(p,w)$ and for large enough $t$, the quantities $\phi_p+w$ and $\phi_p$ are comparable and positive in $\s$. Then one can see that both terms in \eqref{twoterms} are bounded by $C|w|$, uniformly on that neighborhood. Computing on a fixed sufficiently far tail and then on a finite-cylinder shows the continuity of $H_p(w)$.
   	
   	If $(p,w) \in K \times \cX_\delta$ satisfies $\cF(p,w)=0$, $w$ decays in $C^1$ by the $C^1$ weighted embedding, so the associated $\phi$ has the limits $0$ and $\phi_p$, and is nonnegative for sufficiently large positive $t$, by \eqref{eq:phidistance}. Let $\phi_- = \max\{-\phi,0\}$ and take cut-off function $\zeta_R(t) = \zeta(t+R)$ where $R \in \R$, $\zeta$ is smooth, $0 \le \zeta \le 1$, $\zeta(t) = 0$ if $t \le -1$ and $\zeta(t) = 1$ if $t \ge 0$. Testing the equation $\cF(p,w)=0$ with $e^{a_pt}\zeta_R(t)^2\phi_-$ for a fixed $R$ first and then sending $R \to \infty$ yields
   	\[
   	\int_\C e^{a_pt}\big(
   	|\partial_t\phi_-|^2+|\nabla_{\mathbb{S}^{N-1}}\phi_-|^2+b_p\phi_-^2\big)d td\sigma=0.
   	\]
   	Since $N-1+b_p>0$, we get $\phi_-=0$. By the strong maximum principle, we deduce that $\phi>0$, since its positive-end limit is nonzero. The scalar component of \eqref{eq:Fmap} is the phase-normalized condition. Conversely, a normalized connection belongs to the affine weighted space by Lemma~\ref{lem:exptails} and is a zero. Finally, by Proposition~\ref{prop:compact}, we get that the zero set is compact and complete the proof.
   \end{proof}

   Next we study the case when $p \downarrow p_c$. By the angular branch described in \cite[Section~1.2]{DMP}, together with the uniqueness of $\phi_p$,
\begin{equation}\label{eq:phismall}
 \|\phi_p\|_{C^2(\overline \s)}\rightarrow0
 \qquad(p\downarrow p_c).
\end{equation}

\begin{lemma}\label{lem:allsmall}
The supremum of $\|\phi\|_{L^\infty(\C)}$ over all bounded connections at exponent $p$ tends to zero as $p\downarrow p_c$. The same is true for the uniform local $C^1$ norm.
\end{lemma}

\begin{proof}
By Lemma~\ref{lem:bound}, the supremum of $\|\phi\|_{L^\infty(\C)}$ over all bounded connections at exponent $p$ is bounded as $p\downarrow p_c$. Arguing by contradiction, we suppose that there exist connections $\phi_j$, parameters $p_j\downarrow p_c$, and points at which $\phi_j\ge\varepsilon>0$. After translation if necessary, we may assume that $\phi_j(0,z_j)\ge\varepsilon$. Up to a subsequence, we have $\phi_j \to V$ locally in $C^2$, where $V$ is a bounded nonnegative function that is nonzero somewhere. On the other hand, by \eqref{eq:dissipationtotal} and \eqref{eq:phismall},
\[
 \int_\C |(\phi_j)_t|^2=\frac{J_{p_j}}{a_{p_j}}\rightarrow0.
\]
Thus $V$ is stationary and satisfies
\[
 -\Delta_{\mathbb{S}^{N-1}}V-(N-1)V=V^{p_c},\qquad V|_{\partial \s}=0.
\]
Testing with $\varphi_1$ yields $\int_{\s}V^{p_c}\varphi_1=0$, which contradicts $V\ge0$, $V\not\equiv0$. By the uniform local elliptic estimates, we get the same result for the uniform local $C^1$ norm and complete the proof.
\end{proof}

Let
\begin{equation}\label{eq:referenceL}
 L_* =\partial_{tt}+N\partial_t+\Delta_{\cs}+(N-1).
\end{equation}
For $L^2(\s)$-valued functions, we define
\[
 \|h\|_{Z_\eta} :=\sup_{t\in\R}e^{-\eta|t|}
 \big(\|h(t)\|_{L^2}+\|h_t(t)\|_{L^2}\big),
 \quad
 \|R\|_{B_\eta} :=\sup_{t\in\R}e^{-\eta|t|}\|R(t)\|_{L^2}.
\]

\begin{lemma}\label{lem:centerestimate}
Let $h$ be a solution of $L_*h=R$, where $R$ is bounded, $h(t)$ is bounded and axial, has bounded derivative $h_t$ and zero lateral trace. Then, there exist $\eta>0$ and $C_\eta>0$, depending on $N$ and independent of $h$, such that
\begin{equation}\label{eq:centerestimate}
 \| h\|_{Z_\eta}\le C_\eta\big(|\ell(h)|+\| R\|_{B_\eta}\big).
\end{equation}
\end{lemma}

\begin{proof}
Write $h(t)=g(t)\varphi_1+w(t)$, with $w\perp\varphi_1$, and set
\[
 q_1=g',\quad c=g+q_1/N,\quad
 r_1(t)=\int_{\s}R(t,z)\varphi_1(z)d\sigma.
\]
Projecting $L_*h=R$ onto $\text{span} \{\varphi_1\}$ yields $g''+Ng'=r_1$. $g$ needs to be bounded since $h$ is bounded. We solve that
\begin{equation}\label{eq:centerformulas}
 q_1(t)=\int_{-\infty}^t e^{-N(t-s)}r_1(s)d s,
 \qquad
 c(t)=\ell(h)+\frac1N\int_0^t r_1(s) d s.
\end{equation}
On $\text{span}\{\varphi_1\}^{\perp}$, the operator $-\Delta_{\cs}-(N-1)$ is strictly positive. Let $R_\perp$ denote the projection of $R$ onto $\text{span}\{\varphi_1\}^{\perp}$. We have that $L_*w = R_\perp$. Lemma~\ref{lem:green} gives an exponentially decaying kernel $G_\perp$ and its first derivative, with
\[
 w=G_\perp*R_\perp,
 \qquad w_t=G_\perp'*R_\perp.
\]
Choose $\eta$ smaller than $N$ and smaller than the decay rate of these two kernels. Using estimates from the convolution, together with
\[
 \sup_t e^{-\eta|t|}\left|\int_0^t r_1(s)ds \right|
 \le\eta^{-1}\| R\|_{B_\eta},
\]
$g=c-q_1/N$ and $g'=q_1$, we obtain \eqref{eq:centerestimate} and complete the proof.
\end{proof}

\begin{proposition}\label{prop:initialregular}
There exists $\varepsilon_0>0$ such that, for every $p\in(p_c,p_c+\varepsilon_0)$, there exists exactly one normalized axial connection. At this connection, the operator \eqref{eq:phaseindex} is an isomorphism for every sufficiently small admissible positive weight $\delta$.
\end{proposition}

\begin{proof}
After translations, we obtain the existence result by \cite[Proposition~1.1]{DMP}. Moreover, these solutions are axial.
For uniqueness, let $\phi_1,\phi_2$ be two normalized axial connections at the same exponent and let $h=\phi_1-\phi_2$. Then $\ell(h)=0$ and
\begin{equation}\label{eq:differenceatendpoint}
 L_*h=(N-a_p)h_t+ (N-1+b_p) h-(\phi_1^p-\phi_2^p).
\end{equation}
By Lemma~\ref{lem:allsmall}, there exists $M_p \to 0$ such that
\[
 \|(N-a_p)h_t+(N-1+b_p)h-(\phi_1^p-\phi_2^p)\|_{B_\eta}
 \le \big(|N-a_p|+(N-1+b_p)+pM_p^{p-1}\big)\| h\|_{Z_\eta}.
\]
The coefficient at the right-hand side tends to zero as $p\downarrow p_c$. Combining with \eqref{eq:centerestimate} proves $h=0$ for $p$ sufficiently close to $p_c$.

Finally, suppose $h\in\cX_\delta$ satisfies $L_{p,\phi}h=0$ and $\ell(h)=0$. It is bounded with bounded first derivatives, and
\[
 L_*h=(N-a_p)h_t+(N-1+b_p)h-p\phi^{p-1}h.
\]
Reasoning as before, we get $h=0$. By Proposition~\ref{prop:variableindex}, the operator \eqref{eq:phaseindex} has index zero, so it is also surjective. By the bounded inverse theorem, it is indeed an isomorphism. The proof is complete.
\end{proof}

Now we can use a continuation argument to extend the existence result for all $p \in (p_c,p_s)$.

\begin{proposition}\label{prop:existconnection}
For every $p\in(p_c,p_s)$ there exists a bounded normalized axial connection solving \eqref{eq:phi(t,x)}.
\end{proposition}

\begin{proof}
Fix an arbitrary target exponent $P \in (p_c,p_s)$ and choose $p_0<P$ sufficiently close to $p_c$ for Proposition~\ref{prop:initialregular} to apply. Let $K=[p_0,P]$ and we take sufficiently small weight $\delta$ dependent on $K$. By Lemma~\ref{lem:Fproperties}, the map $\cF$ is $C^1$ and its zero set is compact. By Proposition~\ref{prop:variableindex}, its derivative in $w$ is Fredholm of index zero. At $p=p_0$, by Proposition~\ref{prop:initialregular}, the zero is unique and its derivative is an isomorphism. Now we can apply the compact-zero-set continuation lemma, Lemma~\ref{lem:continuation}, to $\cF$ with target space $\cY_\delta\times\R$. This yields a zero at $P$. By Lemma~\ref{lem:Fproperties}, this zero is a normalized connection. Since $P$ was arbitrary, the conclusion follows.
\end{proof}

It remains to show that a connection has a strictly positive, finite Poisson coefficient at its zero end as $t \to -\infty$.

\begin{lemma}\label{lem:preciseasymptotics}
Let $\phi$ be a positive connection at $p\in(p_c,p_s)$. Set
\begin{equation}\label{eq:gammakappa}
 \gamma_p=N-1-\frac{2}{p-1}>0,\qquad
 k_p=\frac{2}{p-1}+1>0.
\end{equation}
Then there exists a constant $C(\phi)>0$ such that
\begin{equation}\label{eq:precisetail}
 \phi(t,z)=C(\phi)e^{\gamma_pt}z_N+o(e^{\gamma_pt})
 \quad\hbox{in }C^1(\overline \s),\quad \text{as } t\to-\infty.
\end{equation}
More precisely,
\begin{equation}\label{eq:coefficient}
 C(\phi)=\frac1{N\|z_N\|_{L^2(\s)}}
 \int_{\R}e^{-\gamma_ps}
 \left(\int_{\s}\phi(s,z)^p\varphi_1(z)d\sigma\right)d s.
\end{equation}
\end{lemma}

\begin{proof}
Write $\gamma=\gamma_p$ and $k=k_p$ for simplicity. Note that
\[
 a_p=k-\gamma,\quad N-1+b_p=\gamma k,\quad \gamma+k=N.
\]
Define
\[
 g(t)=\int_{\s}\phi(t,z)\varphi_1(z)d\sigma,
 \quad F_1(t)=\int_{\s}\phi(t,z)^p\varphi_1(z)d\sigma>0.
\]
Projecting \eqref{eq:phi(t,x)} onto $\text{span}\{\varphi_1\}$ yields
\[
 g''+a_pg'-\gamma k g=-F_1.
\]
Equivalently,
\begin{equation}\label{eq:firstprojectionidentity}
 \left[e^{-\gamma t}(g'+kg)\right]'=-e^{-\gamma t}F_1(t).
\end{equation}
Functions $g$ and $g'$ are bounded as $t \to +\infty$. Thus integration from $t$ to $+\infty$ yields
\begin{equation}\label{eq:firstprojectionformula}
 e^{-\gamma t}(g'(t)+kg(t))
 =\int_t^\infty e^{-\gamma s}F_1(s)d s =: I(t).
\end{equation}
Choose $\beta$ such that $\gamma/p<\beta<\gamma$. By Lemma~\ref{lem:exptails}, we have $\phi=O(e^{\beta t})$ as $t \to -\infty$, so $F_1=O(e^{p\beta t})$. Moreover, $F_1$ is bounded as $t \to +\infty$. It follows that
\begin{equation}\label{eq:Dfinite}
 D:=\int_\R e^{-\gamma s}F_1(s)d s\in(0,\infty).
\end{equation}
 Since $g$ is bounded as $t \to -\infty$,
\[
 g(t)=e^{-kt}\int_{-\infty}^t e^{(k+\gamma)s}I(s)d s.
\]
This, together with $I(t)\to D$ for $t\to-\infty$ and $k+\gamma=N$, proves
\begin{equation}\label{eq:firstmodeleading}
 g(t)=\frac DN e^{\gamma t}+o(e^{\gamma t}), \quad \text{as } t\to-\infty.
\end{equation}

Next, we control the remaining part on $\text{span}\{\varphi_1\}^\perp$. Let $\Pi_\perp$ denote the orthogonal projection onto $\text{span}\{\varphi_1\}^\perp$ and  $w=\Pi_\perp \phi$. It is a bounded solution of
\[
 (\partial_{tt}+a_p\partial_t-B_{-,p})w=-\Pi_\perp(\phi^p).
\]
By this equation, $w$ can be obtained by the convolution of a Green kernel and $-\Pi_\perp(\phi^p)$. Note that the kernel used here is on the subspace $\text{span}\{\varphi_1\}^\perp$. Let $\mu_2$ be the second axial eigenvalue of $B_{-,p}$ and 
\[
\gamma_2 = \frac{-a_p + \sqrt{a_p^2 + 4\mu_2}}{2}.
\]
We have $\mu_2 > b_p + N-1 > 0$, where $b_p + N-1$ is the first eigenvalue of $B_{-,p}$, and $\gamma_2 > \gamma$. The Green kernel decays at least as $e^{\gamma_2t}$ as $t \to -\infty$ and as $e^{-(a_p+\gamma_2)t}$ as $t \to +\infty$. The forcing term is $O(e^{p\beta t})$ in $L^2(\s)$ as $t \to -\infty$ and is bounded as $t \to +\infty$. Splitting the whole-line convolution into $s\le t$, $t<s<0$, and $s\ge0$. For $t<0$, the three terms are bounded by a constant times
\[
 e^{p\beta t},\qquad
 e^{\min\{p\beta,\gamma_2\}t}(1+|t|),\qquad
 e^{\gamma_2t},
\]
respectively. Since $p\beta > \gamma$ and $\gamma_2 > \gamma$, we get $w=o(e^{\gamma t})$ in $L^2(\s)$.
On negative unit slabs, the forcing term is also $O(e^{p\beta t})$ in $L^q$. Recall that $q>N$. By a local estimate like \eqref{eq:localqestimate}, we obtain $w=o(e^{\gamma t})$ in $W^{2,q}$ and hence in $C^1$ by embedding. Finally, combining with \eqref{eq:firstmodeleading} and \eqref{eq:phi1}, we obtain \eqref{eq:precisetail}-\eqref{eq:coefficient}.
\end{proof}

\begin{proof}[Proof of Proposition \ref{prophalfspacesolution}]
Take a connection from Proposition~\ref{prop:existconnection} and set $u(rz)=r^{-\frac{2}{p-1}}\phi(-\log r,z)$. It is positive and solves \eqref{eqhalfspace}. As $r \to 0$, $-\log r \to +\infty$ and $\phi\to\phi_p$ in $C^1(\overline \s)$. A function with zero lateral trace satisfies
\[
 |h(z)|\le\text{dist}(z,\partial \s)\|\nabla_{\cs}h\|_\infty.
\]
Together with \eqref{eq:phidistance}, this turns the $C^1$ convergence into uniform relative convergence and proves \eqref{eq:originratio} for $u$. By Lemma~\ref{lem:preciseasymptotics}, as $r \to +\infty$, $-\log r \to -\infty$ and
\[
 u(rz)=C(\phi)r^{1-N}z_N(1+o(1)).
\]

For $\la > 0$, we define
\[
 u_\lambda(x)=\lambda^{\frac{2}{p-1}}u(\lambda x).
\]
Choose $\lambda=C(\phi)^{1/\gamma_p}$ and call the resulting solution $u_p$. It is straightforward to check that $u_p$ is a solution to problem \eqref{eqhalfspace} satisfying \eqref{eq:originratio} and \eqref{eq:infinityratio}. Finally, by the uniform $C^1$ slab bound for $\phi$ and its zero lateral trace, we get $\phi(t,z)\le C_pz_N$ for all $t$, and hence \eqref{eq:globalbound} holds.
\end{proof}


\begin{remark}[Local very weak trace]
The singularity carries no additional boundary Dirac mass. Indeed, $2/(p-1)<N-1$ implies
\[
 u_p\in L^1_{\mathrm{loc}}(\overline{\R^N_+}),\quad
 x_Nu_p^p\in L^1_{\mathrm{loc}}(\overline{\R^N_+}).
\]
For every $\zeta\in C_c^2(\overline{\R^N_+})$ with $\zeta=0$ on $\partial \R^N_+$, integration by parts on $\R^N_+\setminus\overline B_\varepsilon$ gives, as $\varepsilon\downarrow0$,
\[
 \int_{\R^N_+}\big(u_p\Delta\zeta+u_p^p\zeta\big)d x=0.
\]
In fact, $|\zeta|\le C|x|$ near zero, while $u_p=O(|x|^{-2/(p-1)})$ and $|\nabla u_p|=O(|x|^{-2/(p-1)-1})$ by the cylindrical estimates. The inner hemispherical boundary terms are therefore $O(\varepsilon^{N-1-2/(p-1)})$ and vanish. This is a local very weak statement and no global $L^1(\R^N_+)$ integrability is asserted.
\end{remark}

\begin{remark}
In this section, we do not establish uniqueness for $p$  away from $p_c$ or full non-degeneracy modulo translation. In next section, we develop a different gluing argument from the one in \cite[Section~3.3]{DMP} and do not assume these missing properties.
\end{remark}

\section{The bounded domain case: proof of Theorem \ref{thmmain}}

\subsection{The degree in the full function space with the axial restrictions removed} \label{subdegree}

\begin{lemma}\label{lem:fullangulargap}
On the full space $L^2(\s)$ the operator $B_+=-\Delta_{\cs}+b_p-p\phi_p^{p-1}$ has precisely one negative simple eigenvalue and no zero eigenvalue. On the orthogonal complement of the $O(N-1)$-invariant functions its quadratic form is bounded below by $(a_p+1)\|\cdot\|_{L^2}^2$.
\end{lemma}
\begin{proof}
Write $\alpha=\frac{2}{p-1}$, $a=a_p$, $b=b_p$ for simplicity. For $N\ge3$, we use polar variables
\[
 z=(\sin\theta\,\omega,\cos\theta),\qquad
 \omega\in\mathbb S^{N-2},\quad 0<\theta<\pi/2,
\]
and have
\[
\Delta_{\cs} = \partial_{\theta\theta} + (N-2)\cot \theta \partial_\theta + \frac{1}{\sin^2 \theta}\Delta_{\mathbb{S}^{N-2}}.
\]
Then we use the decomposition in terms of spherical harmonics on $\mathbb S^{N-2}$. Take the standard orthonormal spherical harmonic basis of $L^2(\mathbb{S}^{N-2})$
\[
Y_{l,r}, \quad l = 0,1,\cdots, ~~ 1 \le r \le d_l,
\]
where
\[
-\Delta_{\mathbb{S}^{N-2}}Y_{l,r} = \nu_lY_{l,r}, \quad \nu_l=l(l+N-3).
\]
$l=0$ is the axial sector and $\nu_1=N-2$. In the $l$th sector the operator is
\[
 B_{+,l}=-\partial_{\theta\theta}-(N-2)\cot\theta\,\partial_\theta
 +\frac{\nu_l}{\sin^2\theta}+b-p\phi_p^{p-1}.
\]
Direct computations check that the positive function
\begin{equation}\label{eq:Gangular}
 G(\theta)=\alpha\sin\theta\,\phi_p(\theta)-\cos\theta\,\phi_p'(\theta)
\end{equation}
satisfies
\begin{equation}\label{eq:Gangularidentity}
 B_{+,1}G=(a+1)G.
\end{equation}
Hence, $G(\theta)$ is the principal eigenfunction in the $l=1$ sector. For $l = 0$, by Proposition~\ref{propaxialnondegenerate}, $B_{+,0}$ has exactly one negative simple eigenvalue and no zero eigenvalue. For $l \ge 2$, the first eigenvalue of $B_{+,l}$ is no less than the one of $B_{+,1}$. Thus we can get all the conclusions in this lemma. When $N=2$, note that the nonaxial sector is exactly the odd sector on the semicircle. Then similar to the case when $N \ge 3$, we can arrive at the desired conclusion.
\end{proof}

Fix $q>N$ and let $X_\delta,Y_\delta$ be the spaces in \eqref{eq:Xspace}-\eqref{eq:Yspace} with the axial restrictions removed, that is,
\begin{equation}\label{eq:fullspaces}
 X_\delta=\{w:e^{\delta\rho}w\in W^{2,q}(\C),\ {\rm Tr}_{\R\times\partial \s} w=0\},
 \qquad Y_\delta=\{f:e^{\delta\rho}f\in L^q(\C)\}.
\end{equation}
Set $\phi=\chi\phi_p+w$ and define
\begin{equation}\label{eq:fullmap}
 F_p(w)=\left(\phi_{tt}+a_p\phi_t+\Delta_{\cs}\phi-b_p\phi+(\phi_+)^p,
 \ell(\phi)-\tfrac12m_p\right):X_\delta\rightarrow Y_\delta\times\R.
\end{equation}
The phase $\ell$ and $m_p=\langle\phi_p,\varphi_1\rangle$ are as in \eqref{eq:phase}.

\begin{proposition}\label{prop:fulldegree}
For each fixed $p\in(p_c,p_s)$ and all sufficiently small admissible $\delta>0$, the zero set $K_p=F_p^{-1}(0)$ is compact and the local mod-$2$ Fredholm degree $\deg_2(F_p,K_p)$ is equal to $1$. In particular, for every positive integer $k$,
\[
 \deg_2(F_p^{\times k},K_p^k)=1,
\]
where $F_p^{\times k}$ denotes the product map. The notation denotes the degree in a sufficiently small isolating neighborhood of the indicated compact zero set.
\end{proposition}
\begin{proof}
Choose $p_0\in(p_c,p)$ sufficiently close to $p_c$ and take one small weight $\delta > 0$ dependent on $[p_0,p]$. In view of Lemma~\ref{lem:fullangulargap}, one can check that the nonaxial versions of Proposition~\ref{prop:compact}, Proposition~\ref{prop:variableindex}, and Lemma~\ref{lem:Fproperties} still hold. Moreover, at $p_0$ which is close to $p_c$, the zero of $F_{p_0}$ is unique and regular in the full space as well. Thus the full-space map at $p_0$ has degree $1$. Applying the finite-dimensional reduction in Appendix~\ref{sec:continuationlemma} yields $\deg_2(F_p,K_p) = 1$. Finally, for a positive integer $k$, $\deg_2(F_p^{\times k},K_p^k)=1$ follows by taking products of the finite-dimensional reductions.
\end{proof}

\begin{lemma}[Persistence under small full-space perturbations]\label{lem:persistence}
Let $F:X\to Y$ be $C^1$, Fredholm of index zero near a compact zero set $K$, and suppose that its local mod-$2$ degree at $K$ is nonzero. Then, there exists a bounded neighborhood $U$ of $K$ such that, if $R_\varepsilon:U\to Y$ is $C^1$ and
\[
 \|R_\varepsilon\|_{C^1(U)}\rightarrow0,
\]
then $F+R_\varepsilon$ has a zero $u_\ep$ in $U$ for every sufficiently small $\varepsilon$ and up to a subsequence, $u_\ep \to u$ in $X$ as  $\varepsilon\to0$ for some $u \in K$.
\end{lemma}
\begin{proof}
Similar to the proof of Lemma~\ref{lem:continuation}, we can take a finite-dimensional $E\subset Y$ such that $\text{Ran } DF(x)+E=Y$ near $K$. Let $\pi:Y\to Y/E$ be the quotient map. By implicit function theorem for the map $\pi F$, the set $M=\{\pi F=0\}$ is a finite-dimensional $C^1$ manifold near $K$, of dimension $\dim E$. On it $f=F|_M$ takes values in $E$. We choose a relatively compact isolating neighborhood $D$ of $K$ in $M$ whose boundary has positive distance from zero under $f$.

On \(\overline D\), choose continuous, uniformly bounded right inverses \(A_x:Y/E\to X\) of \(\pi DF(x)\). Using a uniform contraction argument to solve
\[
\pi(F+R_\varepsilon)(x+A_xy)=0,
\]
we obtain a continuous \(y=y_\varepsilon(x)\), with
\(\sup_{\overline D}\|y_\varepsilon(x)\|\to0\).
As $\ep \to 0$, the reduced maps
\[
f_\varepsilon(x)
=(F+R_\varepsilon)(x+A_xy_\varepsilon(x))\in E
\]
converge uniformly to \(f\). Since the finite-dimensional mod-2 degree of $f$ on $D$ is nonzero and the perturbation $\|R_\varepsilon\|_{C^1(U)}$ is small, one can check that the degree of $f_\ep$ on $D$ is also nonzero when $\ep$ is small enough, so there exists a zero of $f_\ep$. This gives the solution $u_\ep$ of $F+R_\varepsilon$. Finally, using the compactness of \(\overline D\) and the uniform convergence \(y_\varepsilon\to0\), we obtain the convergence of $u_\ep$ up to a subsequence.
\end{proof}

\subsection{An exact inner-outer system on a bounded domain} \label{subinner-outer}

We now focus on the proof of Theorem~\ref{thmmain}. Fix $p\in(p_c,p_s)$ and abbreviate $\alpha=2/(p-1)$, $a=a_p$, $b=b_p$. Let $y_1,\ldots,y_k$ be distinct points of the boundary of a bounded connected $C^{2,\beta}$ domain $\Omega$. Fix $q>N$ as above. All constants below may depend on $p,\Omega,\{y_i\}$, the weight $\delta$, and a bounded neighborhood of $K_p^k$, but not on the small scale $\varepsilon$.

For each $i$, choose a boundary-flattening chart $\Phi_i$ on a fixed half-ball, with $\Phi_i(0)=y_i$, $D\Phi_i(0)$ orthogonal, and the last coordinate pointing inward. Such charts can be chosen as a rotation followed by the graph map $(\xi',\xi_N)\mapsto(\xi',\xi_N+f_i(\xi'))$, where $f_i(0)= 0$, $\nabla f_i(0)=0$. In these coordinates, following directly by differentiating the graph map we obtain
\begin{equation}\label{eq:metricerror}
 (\Delta_x B)\circ\Phi_i=\Delta_\xi(B\circ\Phi_i)
 + \sum_{l,m}a_i^{lm}(\xi)\partial_{lm}(B\circ\Phi_i)
 + \sum_{l}c_i^l(\xi)\partial_l(B\circ\Phi_i),
 \qquad |a_i^{lm}(\xi)|\le C|\xi|,\quad |c_i^l(\xi)|\le C.
\end{equation}

By decreasing the admissible weight in Proposition~\ref{prop:fulldegree} if necessary, we can assume that
\begin{equation}\label{eq:extradelta}
 0<\delta<\min\left\{1,\ N-1-\alpha,\ \frac{\alpha+1}{p+1}\right\}.
\end{equation}
Take
\begin{equation}\label{eq:kappachoice}
 0<\kappa<\min\left\{
 1,\ \frac{\delta}{\alpha+\delta+2},\
 \frac{(p-1)\delta}{p+(p-1)\delta}\right\},
 \qquad R=\varepsilon^\kappa.
\end{equation}
Then $\varepsilon/R\to0$ and $R\to0$.

For $i$ we take a smooth cut-off function $\eta_i(x)$, which equals $1$ for $|\Phi_i^{-1}(x)|\le R$, vanishes for $|\Phi_i^{-1}(x)|\ge2R$ and outside the chart, and has the usual derivative bounds $C_jR^{-j}$ for $j\le2$. All supports with different $i$ are disjoint when $\varepsilon$ is small.  Moreover, we will use a second smooth cut-off function $\zeta_i(\xi)$, equal to $1$ for $|\xi|\le2R$ and zero for $|\xi|\ge3R$.

Take a bounded neighborhood $U$ of $K_p^k$ in $X_\delta^k$. For $\boldsymbol w=(w_1,\ldots,w_k) \in U$ we set
\begin{equation}\label{eq:rawprofiles}
 \phi_i=\chi\phi_p+w_i,\quad
 B_i(\Phi_i(rz))=r^{-\alpha}\phi_i\big(\log(\varepsilon/r),z\big),
 \quad r>0,\ z\in \s,
\end{equation}
where the functions $B_i$ are used only inside their charts. Note that the functions $\phi_i$ need not be axially symmetric. We would like to search for solutions in the following form:
\begin{equation}\label{eq:exactansatz}
 u=\sum_{i=1}^k\eta_i B_i+h,
 \qquad h\in Z:=W^{2,q}(\Omega)\cap W^{1,q}_0(\Omega).
\end{equation}
Let 
\begin{equation}\label{eq:commutator}
 \varpi_i=2\nabla\eta_i\cdot\nabla B_i+B_i\Delta\eta_i,
\end{equation}
extended by zero away from \{$R<|\Phi_i^{-1}(x)|<2R$\}. The outer error is given by the equation
\begin{equation}\label{eq:outereq}
 -\Delta h =
 \left(\sum_i\eta_iB_i+h\right)_+^p
 -\sum_i\eta_i(B_i+h)_+^p+\sum_i \varpi_i,
 \qquad h|_{\partial\Omega}=0.
\end{equation}
For the inner equation, set $t=\log(\varepsilon/r)$ and define
\begin{equation}\label{eq:Hdefinition}
 H_i(t,z)=r^\alpha h(\Phi_i(rz)),\qquad r=\varepsilon e^{-t},
\end{equation}
where it is needed when $r<3R$. The geometric error operator $\mathcal E_{i,\varepsilon}$ is defined by
\begin{equation}\label{eq:innergeometry}
 r^{\alpha+2}(\Delta_x B_i)(\Phi_i(rz))
 =\mathcal L_0\phi_i+\mathcal E_{i,\varepsilon}\phi_i,
 \quad \mathcal L_0=\partial_{tt}+a\partial_t+\Delta_{\cs}-b.
\end{equation}
In smooth angular coordinate patches, $\mathcal E_{i,\varepsilon}$ is an operator of order at most two with coefficients bounded by $Cr$. The inner system is
\begin{equation}\label{eq:innereq}
 \begin{split}
 \mathcal L_0\phi_i+(\phi_i)_+^p+
 \zeta_i(rz)[\mathcal E_{i,\varepsilon}\phi_i+(\phi_i+H_i)_+^p-(\phi_i)_+^p]&=0,\\
 \ell(\phi_i)&=\tfrac12m_p,
 \qquad i=1,\ldots,k.
 \end{split}
\end{equation}
The cut-off functions $\zeta_i$ are extended by zero where the charts are not used.

\begin{lemma}[Exactness of the coupled equations]\label{lem:exactness}
A solution of \eqref{eq:outereq} and \eqref{eq:innereq} gives, by \eqref{eq:exactansatz}, a distributional solution of
\[
 \Delta u+(u_+)^p=0\quad\hbox{in }\Omega,
\]
which is locally classical away from the prescribed boundary points and has zero boundary values there.
\end{lemma}
\begin{proof}
On the support of $\eta_i$, $\zeta_i = 1$. Since $\alpha p=\alpha+2$, the first equation in \eqref{eq:innereq} and \eqref{eq:innergeometry} yield
\[
 \Delta_x B_i=-(B_i+h)_+^p.
\]
Hence, by \eqref{eq:outereq},
\[
 \Delta\left(\sum_i\eta_iB_i+h\right)
 =-\sum_i\eta_i(B_i+h)_+^p+\sum_i \varpi_i+\Delta h
 =-\left(\sum_i\eta_iB_i+h\right)_+^p.
\]
On each compact subset away from the boundary points, all terms are locally $W^{2,q}$. By a standard bootstrap argument, we get that the solution is locally classical. The zero lateral trace of $\phi_i$, and the trace of $h$, prove the boundary condition.
\end{proof}

Set
\begin{equation}\label{eq:smallm}
 M_\varepsilon=\varepsilon^\delta R^{-\alpha-\delta},
 \quad m_\varepsilon=\varepsilon^\delta R^{-\alpha-\delta-2}.
\end{equation}
By \eqref{eq:kappachoice}, $m_\varepsilon\to0$. Uniformly on bounded subsets of $X_\delta$, the weighted Sobolev embedding yields
\begin{equation}\label{eq:vpointwise}
 |\phi_i(t,z)|+|\nabla_{t,\cs}\phi_i(t,z)|
 \le C e^{\delta t}~~(t\le-1),
 \quad |\phi_i|+|\nabla_{t,\cs}\phi_i|\le C~~(t\ge-1).
\end{equation}
Consequently, when $R<r<2R$ we have
\begin{equation}\label{eq:annularbound}
 |B_i|+R|\nabla_xB_i|\le C M_\varepsilon,
 \qquad \|\sum_i \varpi_i\|_{L^q(\Omega)}\le C m_\varepsilon.
\end{equation}
For any $\boldsymbol w \in U$, $\boldsymbol w' \in X_\delta^k$, denote 
\begin{align} \label{eqdotphi}
	\dot{\phi}_i = D_{w_i}\phi_i(w_i)[w_i'], \quad \dot{B}_i = D_{w_i}B_i(w_i)[w_i'], \quad \dot{\varpi}_i = D_{w_i}\varpi_i(w_i)[w_i'].
\end{align}
We also have 
\begin{equation}\label{eq:vpointwisederivative}
 |\dot{\phi}_i|+|\nabla_{t,\cs}\dot{\phi}_i|
 \le C e^{\delta t}\|w_i'\|_{X_\delta}~~(t\le-1),
 \quad |\dot{\phi}_i|+|\nabla_{t,\cs}\dot{\phi}_i|\le C\|w_i'\|_{X_\delta}~~(t\ge-1),
\end{equation}
and when $R<r<2R$,
\begin{equation}\label{eq:annularboundderivative}
 |\dot{B}_i|+R|\nabla_x\dot{B}_i|\le C M_\varepsilon\|w_i'\|_{X_\delta},
 \quad \|\sum_i \dot{\varpi}_i\|_{L^q(\Omega)}\le C m_\varepsilon\|\boldsymbol w'\|_{X_\delta^k}.
\end{equation}

\begin{proposition}\label{prop:outer}
For every sufficiently small $\varepsilon$ and each $\boldsymbol w\in U$, equation~\eqref{eq:outereq} has a unique solution, denoted by $h_\varepsilon(\boldsymbol w)$, in the ball
\[
 \{\|h\|_Z\le C m_\varepsilon\}.
\]
The map $\boldsymbol w \mapsto h_\varepsilon(\boldsymbol w)$ is $C^1$ and satisfies
\begin{equation}\label{eq:outerC1}
 \sup_{\boldsymbol w\in U}\left(
 \|h_\varepsilon(\boldsymbol w)\|_Z+
 \|D h_\varepsilon(\boldsymbol w)\|_{\mathcal L(X_\delta^k, Z)}
 \right)\le C m_\varepsilon.
\end{equation}
\end{proposition}
\begin{proof}
The Dirichlet inverse $T=(-\Delta_\Omega)^{-1}:L^q(\Omega)\to Z$ is bounded since $\Omega$ is bounded, see, for example, \cite[Chapter~9]{GT}. Write the right side of \eqref{eq:outereq} as $N(\boldsymbol w,h)+\sum_i \varpi_i$, that is, 
\[
N(\boldsymbol w,h) = \left(\sum_i\eta_iB_i+h\right)_+^p
 -\sum_i\eta_i(B_i+h)_+^p.
\]
Note that $N=0$ for $|\Phi_i^{-1}(x)|\le R$, $B_i$ are bounded by $CM_\varepsilon$ for $R < |\Phi_i^{-1}(x)| < 2R$, and away from the cut-off supports $N=(h_+)^p$. Therefore, by $Z\hookrightarrow C^1(\overline\Omega)$ and the inequality
$|a_+^p-b_+^p|\le C_p(|a|+|b|)^{p-1}|a-b|$, 
\begin{align*}
 \|N(\boldsymbol w,h)\|_{L^q}&\le C(M_\varepsilon^p+\|h\|_Z^p),\\
 \|D_hN(\boldsymbol w,h)\|_{\mathcal L(Z,L^q)}
 &\le C(M_\varepsilon+\|h\|_Z)^{p-1}.
\end{align*}
Since $M_\varepsilon=R^2m_\varepsilon$ and $m_\varepsilon\to0$, the map $h\mapsto T(N+\sum_i\varpi_i)$ maps $\{\|h\|_Z\le C m_\varepsilon\}$ into itself for sufficiently small $\varepsilon$ and is a contraction there, uniformly on $U$. This solves $h_\varepsilon(\boldsymbol w)$ in $\{\|h\|_Z\le C m_\varepsilon\}$.

For fixed $\varepsilon$, considering for $|\Phi_i^{-1}(x)|\le R$, $R \le |\Phi_i^{-1}(x)|\le 2R$ and outside the supports of $\eta_i$ respectively, one checks that $N$ is a $C^1$ map into $L^q$. By \eqref{eq:annularbound}, $\sum_i\varpi_i$ is also a $C^1$ map into $L^q$. In view of \eqref{eq:annularboundderivative}, their Fr\'echet derivatives with respect to $\boldsymbol w$ have norm at most $Cm_\varepsilon$. Differentiating the fixed-point equation $h = T(N+\sum_i\varpi_i)$ yields
\[
 D h_\varepsilon=(I-TD_hN)^{-1}T
 \left(D_{\boldsymbol w}N+D_{\boldsymbol w}\sum_i \varpi_i\right).
\]
It is straightforward to check that the inverse $(I-TD_hN)^{-1}$ has a uniform bound on $U$, proving \eqref{eq:outerC1}.
\end{proof}

Let 
\[
\mathfrak R_\varepsilon = \left((\mathfrak R_{1,\varepsilon}(\boldsymbol w),0), \cdots, (\mathfrak R_{k,\varepsilon}(\boldsymbol w),0)\right) : U\rightarrow(Y_\delta\times\R)^k,
\]
where
\[
\mathfrak R_{i,\varepsilon}(\boldsymbol w) = \zeta_i(rz)[\mathcal E_{i,\varepsilon}\phi_i+(\phi_i+H_i)_+^p-(\phi_i)_+^p].
\]
Inserting $h=h_\varepsilon(\boldsymbol w)$ into \eqref{eq:innereq} enables us to introduce the following map
\begin{equation}\label{eq:reducedmap}
 \mathfrak F_\varepsilon(\boldsymbol w)
 =F_p^{\times k}(\boldsymbol w)+\mathfrak R_\varepsilon(\boldsymbol w)
 :U\rightarrow(Y_\delta\times\R)^k.
\end{equation}

\begin{proposition}\label{prop:innersmall}
With \eqref{eq:extradelta}-\eqref{eq:kappachoice}, as $\varepsilon \to 0$,
\begin{equation}\label{eq:C1small}
 \|\mathfrak R_\varepsilon\|_{C^1(U,(Y_\delta\times\R)^k)}\rightarrow0.
\end{equation}
\end{proposition}
\begin{proof}
\textit{Step 1: the geometric error.} On the support of $\zeta_i$ the coefficients of $\mathcal E_{i,\varepsilon}$ are bounded by $CR$, and it follows that
\[
 \|\zeta_i\mathcal E_{i,\varepsilon}w_i\|_{Y_\delta}
 \le CR\|w_i\|_{X_\delta}.
\]
Recall that $r=\varepsilon e^{-t}$. The term $\chi\phi_p$ appears only when $t\ge-1$, and the $L^q(\s)$ norm of $\zeta_i\mathcal E_{i,\varepsilon}\chi\phi_p$ is bounded by $C\varepsilon e^{-t}$. Since $\delta<1$, a direct computation yields $\|\zeta_i\mathcal E_{i,\varepsilon}\chi\phi_p\|_{Y_\delta}=O(\varepsilon)$. Hence, 
\[
\|\zeta_i\mathcal E_{i,\varepsilon}\phi_i\|_{Y_\delta} \le C\varepsilon + CR\|w_i\|_{X_\delta}, 
\]
and so
\[
\sup_{\boldsymbol w \in U}\|\zeta_i\mathcal E_{i,\varepsilon}\phi_i\|_{Y_\delta} \le C_U(R+\varepsilon).
\]
Define a map
\[
\mathcal{G}_{i,\varepsilon}: \boldsymbol w \mapsto \zeta_i\mathcal E_{i,\varepsilon}\phi_i.
\]
Since
\[
D\mathcal{G}_{i,\varepsilon}(\boldsymbol w)[\boldsymbol w'] = \zeta_i\mathcal E_{i,\varepsilon} w_i',
\]
we have
\[
\|D\mathcal{G}_{i,\varepsilon}(\boldsymbol w)\|_{\mathcal L(X_\delta^k,Y_\delta)} \le CR.
\]
Combining above arguments, we get
\[
\|\mathcal{G}_{i,\varepsilon}\|_{C^1(U,Y_\delta)} \le C_U(R+\varepsilon) \to 0 \quad \text{as } \varepsilon \to 0.
\]

\textit{Step 2: size of the lifted outer correction.} Using \eqref{eq:outerC1} and
\[
|h(\Phi_i(rz))|\le C\| h\|_Z\text{dist}(\Phi_i(rz),\partial\Omega)\le C\| h\|_Z r,
\]
we deduce
\begin{equation}\label{eq:Hbound}
 |H_i(t,z)|\le C m_\varepsilon r^{\alpha+1},
 \quad
 |D_{\boldsymbol w}H_i[\boldsymbol w']|
 \le C m_\varepsilon r^{\alpha+1}\|\boldsymbol w'\|_{X_\delta^k},
 \quad r=\varepsilon e^{-t}<3R.
\end{equation}
Let $n(s)=(s_+)^p$. We have
\begin{equation}\label{eq:Ndiffbound}
 |n(\phi_i+H_i)-n(\phi_i)|\le C(|\phi_i|^{p-1}|H_i|+|H_i|^p).
\end{equation}
When $t \ge 0$, $\phi_i$ is uniformly bounded. Since $\delta<\alpha+1$, by \eqref{eq:Hbound} and direct computations, we obtain
\[
 \|\zeta_i[n(\phi_i+H_i)-n(\phi_i)]\|_{Y_\delta(\{t\ge0\})}
 \le C m_\varepsilon\varepsilon^{\alpha+1}
\]
for small $\varepsilon$.

Recall that $\zeta_i(rz) = 0$ if $r \ge 3R$. Hence, for $t < 0$, it suffices to consider $-L_\varepsilon<t<0$, where $-L_\varepsilon=\log(\varepsilon/(3R))$. By \eqref{eq:vpointwise} and \eqref{eq:Hbound}, we derive that
\begin{align*}
 e^{-\delta t}|\phi_i|^{p-1}|H_i|
 &\le C m_\varepsilon\varepsilon^{\alpha+1}
 e^{-(\alpha+1-(p-2)\delta)t},\\
 e^{-\delta t}|H_i|^p
 &\le C m_\varepsilon^p\varepsilon^{p(\alpha+1)}
 e^{-(p(\alpha+1)+\delta)t}.
\end{align*}
By \eqref{eq:extradelta}, $\alpha+1-(p-2)\delta > 0$, $p(\alpha+1)+\delta > 0$. Integrating their $q$th powers over $(-L_\varepsilon,0)\times \s$ yields
\begin{align}
 \|\zeta_i[n(\phi_i+H_i)-n(\phi_i)]\|_{Y_\delta}
 \le C\big[&m_\varepsilon\varepsilon^{\alpha+1}
 +m_\varepsilon\varepsilon^{(p-2)\delta}R^{\alpha+1-(p-2)\delta}\notag\\
 &+m_\varepsilon^p\varepsilon^{-\delta}R^{p(\alpha+1)+\delta}\big].
 \label{eq:innerintegrated}
\end{align}
By \eqref{eq:smallm},
\[
 m_\varepsilon\varepsilon^{(p-2)\delta}R^{\alpha+1-(p-2)\delta}=\varepsilon^{(p-1)\delta}R^{-1-(p-1)\delta},\quad
 m_\varepsilon^p\varepsilon^{-\delta}R^{p(\alpha+1)+\delta} = \varepsilon^{(p-1)\delta}R^{-p-(p-1)\delta}.
\]
As $p>1$ and $R<1$, by \eqref{eq:kappachoice}, the right side in \eqref{eq:innerintegrated} tends to $0$ as $\varepsilon \to 0$.

\textit{Step 3: derivatives of the nonlinear error.} Define a map
\[
\mathcal{N}_{i,\varepsilon}: \boldsymbol w \mapsto \zeta_i[n(\phi_i+H_i)-n(\phi_i)].
\]
It remains to estimate the Fr\'echet derivative
\begin{equation}\label{eq:derivativeerror}
 D\mathcal{N}_{i,\varepsilon}(\boldsymbol w)[\boldsymbol w'] = \zeta_i
 [n'(\phi_i+H_i)-n'(\phi_i)]w_i' + \zeta_i n'(\phi_i+H_i)D H_i[\boldsymbol w'].
\end{equation}
Let $\tau=\min\{1,p-1\}$. Since $\phi_i$ is uniformly bounded and
\[
 \sup_{\boldsymbol w \in U, r < 3R}|H_i|\le C m_\varepsilon R^{\alpha+1}\rightarrow0,
\]
the scalar $C^{0,\tau}$ estimate for $n'$ gives
\[
 \sup_{\boldsymbol w \in U, r < 3R}|n'(\phi_i+H_i)-n'(\phi_i)|
 \le C(m_\varepsilon R^{\alpha+1})^\tau=o(1).
\]
Thus we can obtain that
\[
\sup_{\boldsymbol w \in U}\|\zeta_i
 [n'(\phi_i+H_i)-n'(\phi_i)]\|_{\mathcal L(X_\delta^k,Y_\delta)}\le C(m_\varepsilon R^{\alpha+1})^\tau=o(1).
\]
Moreover, using \eqref{eq:Hbound} we get
\[
|\zeta_i n'(\phi_i+H_i)D H_i[\boldsymbol w']| \le C\zeta_i(|\phi_i|^{p-1}m_\varepsilon r^{\alpha+1}+m_\varepsilon^p r^{p(\alpha+1)})\|\boldsymbol w'\|_{X_\delta^k}.
\]
As in Step 2, integrating on $\{t \ge 0\}$ and $\{-L_\varepsilon < t < 0\}$ respectively yields
\begin{align}
	\|\zeta_i n'(\phi_i+H_i)D H_i[\boldsymbol w']\|_{Y_\delta} \le C\big[&m_\varepsilon\varepsilon^{\alpha+1}
 +m_\varepsilon\varepsilon^{(p-2)\delta}R^{\alpha+1-(p-2)\delta}\notag\\
 &+m_\varepsilon^p\varepsilon^{-\delta}R^{p(\alpha+1)+\delta}\big]\|\boldsymbol w'\|_{X_\delta^k}.
\end{align}
Combining above arguments yields
\[
\sup_{\boldsymbol w \in U} \|D\mathcal{N}_{i,\varepsilon}(\boldsymbol w)\|_{\mathcal L(X_\delta^k,Y_\delta)} \to 0 \quad \text{as } \varepsilon \to 0.
\]
Together with Step 2, we get 
\[
\|\mathcal{N}_{i,\varepsilon}\|_{C^1(U,Y_\delta)} \to 0 \quad \text{as } \varepsilon \to 0.
\]
Hence, \eqref{eq:C1small} holds.
\end{proof}

\begin{proposition}\label{prop:gluedsolutions}
For every sufficiently small $\varepsilon>0$, there are profiles $\boldsymbol w_\varepsilon\in U$ being a zero point of \eqref{eq:reducedmap}, and a corresponding solution $u_\varepsilon$ of \eqref{eq:exactansatz}. The profiles can be chosen precompact in $X_\delta^k$ as $\varepsilon\to0$, with all limits in $K_p^k$, and
\[
 \|h_\varepsilon\|_Z\le C m_\varepsilon\rightarrow0.
\]
\end{proposition}
\begin{proof}
By Proposition~\ref{prop:fulldegree}, $F_p^{\times k}$ has local degree one at the compact zero set $K_p^k$. Choose its persistence neighborhood $U$ from Lemma~\ref{lem:persistence}. In view of Proposition~\ref{prop:innersmall}, Lemma~\ref{lem:persistence}, Proposition~\ref{prop:outer}, Lemma~\ref{lem:exactness}, we can complete the proof.
\end{proof}

\subsection{Proof of Theorem~\ref{thmmain}} \label{subproof2}

\begin{proof}[Proof of Theorem~\ref{thmmain}]
For $p\in(p_c,p_s)$ take the solutions from Proposition~\ref{prop:gluedsolutions}. We check that these solutions satisfy the following properties.

\textit{Strong singularity and positivity.} In a sufficiently small core of $y_i$, the cut-off function $\eta_i = 1$ and
\[
 r^\alpha u_\varepsilon(\Phi_i(rz))
 =\phi_{i,\varepsilon}(\log(\varepsilon/r),z)
 +r^\alpha h_\varepsilon(\Phi_i(rz)).
\]
Note that when $t \ge 1$,
\[
 |\phi_{i,\varepsilon}(t,z)-\phi_p(z)| = |w_{i,\varepsilon}| \le C'e^{-\delta t}\text{dist }(z,\partial \s) \le Ce^{-\delta t}z_N.
\]
We also have $|h_\varepsilon(\Phi_i(rz))|\le Cm_\varepsilon r z_N$. Since $C_1 z_N \le \phi_p(z) \le C_2 z_N$, for $r<\varepsilon/e$ we obtain
\begin{equation}\label{eq:relativebounded}
 \sup_{z\in \s}\left|
 \frac{r^\alpha u_\varepsilon(\Phi_i(rz))}{\phi_p(z)}-1
 \right|
 \le C(r/\varepsilon)^\delta+Cm_\varepsilon r^{\alpha+1}.
\end{equation}
Thus $u_\varepsilon$ is positive near each prescribed point $y_i$. On the domain with small half-neighborhoods of these points removed, $-\Delta u_\varepsilon=(u_\varepsilon)_+^p\ge0$, the remaining boundary values are zero, and the artificial inner boundary values are nonnegative. By the maximum principle, $u_\varepsilon\ge0$. Then using the strong maximum principle we obtain that $u_\varepsilon$ is positive and the equation of $u_\varepsilon$ becomes $-\Delta u_\varepsilon=u_\varepsilon^p$. Local boundary regularity shows that $u_\varepsilon \in C^2_{\mathrm{loc}}(\overline\Omega\setminus\{y_i\})$. On every fixed nontangential cone, $\phi_p$ has a positive lower bound, so $u_\varepsilon(x) \to +\infty$ as $x \to y_i$ nontangentially.

\textit{The very weak formulation.} For fixed $\varepsilon$, near $y_i$ the previous estimates and local elliptic regularity yield
\begin{equation}\label{eq:physicalgrowth}
 u_\varepsilon(x)\le C_\varepsilon|x-y_i|^{-\alpha},
 \qquad |\nabla u_\varepsilon(x)|\le C_\varepsilon|x-y_i|^{-\alpha-1}.
\end{equation}
Since $\alpha<N-1$, both $u_\varepsilon$ and $d(x)u_\varepsilon^p$ belong to $L^1(\Omega)$, where $d(x)=\text{dist}(x,\partial\Omega)$. For $f\in C^2(\overline\Omega)$ with $f = 0$ on $\partial \Omega$, $|f|\le Cs$ on a small artificial inner boundary of radius $s$. Integrating the classical equation on the punctured domain and then sending $s\downarrow0$, we get
\[
 \int_\Omega(u_\varepsilon\Delta f +u_\varepsilon^p f) d x=0.
\]
This proves the zero very weak trace.

\textit{Sobolev trace and concentration.} Since profiles $\boldsymbol w_\varepsilon$ are in a bounded subset of $X_\delta^k$, when $r<2R$, the embedding and a direct computation yield
\[
 |B_{i,\varepsilon}|\le Cr^{-\alpha},\qquad
 |\nabla(\eta_iB_{i,\varepsilon})|\le Cr^{-\alpha-1},
\]
with a constant independent of $\varepsilon$. Therefore, for every
$1<s<N/(\alpha+1)$,
\begin{equation}\label{eq:Sobolevsmall}
 \|\sum_i\eta_iB_{i,\varepsilon}\|_{W^{1,s}(\Omega)}
 \le C_s R^{N/s-\alpha-1}\rightarrow0 \quad \text{as } \varepsilon \to 0.
\end{equation}
Next we check $\sum_i\eta_iB_{i,\varepsilon} \in W_0^{1,s}(\Omega)$. For this, truncate the profiles with smooth cut-off functions $\varsigma_i$ so that $\varsigma_i\eta_iB_{i,\varepsilon} = 0$ when $r \le \nu$ and $\varsigma_i\eta_iB_{i,\varepsilon} = \eta_iB_{i,\varepsilon}$ when $r \ge 2\nu$, with a radius $\nu > 0$. Note that $\sum_i\varsigma_i\eta_iB_{i,\varepsilon} \in W_0^{1,s}(\Omega)$. Letting $\nu \to 0$ and computing directly yields that $\eta_iB_{i,\varepsilon} \in W^{1,s}_0(\Omega)$. Recall that the outer correction $h_\varepsilon$ tends to zero in $W^{2,q}$, and hence in $W^{1,s}$. Since $s < N/(\alpha+1) < N <q$ and $\Omega$ is bounded, we obtain 
\[
u_\varepsilon \in W_0^{1,s}(\Omega), \quad \|u_\varepsilon\|_{W^{1,s}(\Omega)} \to 0 ~~ \text{as } \varepsilon \to 0.
\]
Denote $B^+_{2R} = \{\xi: |\xi|<2R, \xi_N > 0\}$. It is straightforward to check that
\[
 \int_{\Omega\cap\bigcup_i\Phi_i(B^+_{2R})}d(x)u_\varepsilon^p d x
 \le C (R^{N-1-\alpha}+ m_\varepsilon^pR^{N+1}).
\]
Moreover, outside $\bigcup_i\Phi_i(B^+_{2R})$, $\eta_iB_i = 0$ and $u_\varepsilon = h_\varepsilon$, implying
\[
\int_{\Omega\setminus (\bigcup_i\Phi_i(B^+_{2R}))}d(x)u_\varepsilon^p d x
 \le Cm_\varepsilon^p.
\]
Since $u_\varepsilon$ is positive, we conclude
\[
\|d(x)u_\varepsilon^p\|_{L^1(\Omega)} \to 0 \quad \text{as } \varepsilon \to 0.
\]

On every compact subset of $\overline\Omega\setminus\{y_i\}$, the cut-off terms vanish for all sufficiently small $\varepsilon$, so $u_\varepsilon=h_\varepsilon\to0$ in $C^1$. The local boundary Schauder estimates improve this to $C^2$ convergence. Solutions $u_\varepsilon$ are nonzero, but converge locally to zero. Hence, there are infinitely many distinct solutions.

Finally, on smooth domains, the case when $p=p_c$ has been solved in \cite{DMP}, see Theorem~\ref{thm:dmp}. The proof is complete.
\end{proof}

\appendix
\section{Continuation with a compact zero set}\label{sec:continuationlemma}

We give the topological principle used in Proposition~\ref{prop:existconnection}. The formulation is local near the zero set and does not assume that the nonlinear map is proper on an arbitrary infinite-dimensional ball. This is the local mod-$2$ continuation argument underlying the $C^1$ Fredholm degree of \cite{PR}.

\begin{lemma}\label{lem:continuation}
Let $X,Y$ be real Banach spaces and let $F:[a,b]\times X\to Y$ be $C^1$, with a $C^1$ extension to a neighborhood of $[a,b]$. Suppose that $D_xF(s,x)$ is Fredholm of index zero in a neighborhood of the zero set
\[
 Z=\{(s,x)\in[a,b]\times X:F(s,x)=0\},
\]
and $Z$ is compact. Assume that $F(a,\cdot)$ has exactly one zero $x_a$, and $D_xF(a,x_a)$ is an isomorphism. Then $F(s,\cdot)$ has a zero for every $s\in[a,b]$.
\end{lemma}

\begin{proof}
It is sufficient to prove the result for $b$. Suppose by contradiction that $F(b,\cdot)$ has no zero.

Fix $z_0 = (s_0,x_0) \in Z$, we can choose a finite-dimensional subspace $E_{z_0} \subset Y$ such that 
\[
Y = \text{Ran}D_xF(s_0,x_0) \oplus E_{z_0}.
\]
Let $\pi_{z_0} : Y \to Y/E_{z_0}$ be the quotient map. Then $T_0 = \pi_{z_0}D_xF(s_0,x_0)$ is surjective and has a finite-dimensional kernel. Hence, we can take a bounded right inverse $A_0 : Y/E_{z_0} \to X$ such that $T_0A_0 = I$. For $z = (s,x)$ near $z_0$, denote $T_z = \pi_{z_0}D_xF(s,x)$. Since $F$ is $C^1$, $T_zA_0$ is invertible and $A_0(T_zA_0)^{-1}$ is the right inverse of $T_z$. Therefore, $T_z$ is also surjective, that is,
\[
\text{Ran}D_xF(s,x) + E_{z_0} = Y.
\]
Since $Z$ is compact, we can take the sum of finitely many $E_{z_j}$, $j = 1,\cdots,m$, denoted by $E\subset Y$, such that
\[
 \text{Ran } D_xF(s,x)+E=Y
\]
throughout an open neighborhood $U$ of $Z$. Clearly, $E$ is finite-dimensional. Let $\pi:Y\to Y/E$ be the quotient map. Applying the implicit function theorem, we get that the equation $\pi F=0$ defines in $U$ a finite-dimensional $C^1$ manifold
\[
 M=\{(s,x)\in U:\pi F(s,x)=0\},
 \qquad \dim M=\dim E+1,
\]
where we have used that $D_xF(s,x)$ is Fredholm of index zero in a neighborhood of. Let $P : M \to \R$ be defined by $P(s,x) = s$. Note that $\pi D_xF$ is onto. For each $\sigma \in \R$, we can find $h \in X$ such that $(\sigma,h) \in T_{s,x}M$, that is, $DP$ is a surjective. For each $c \in [a,b]$, $M_c = M \cap (\{c\}\times X)$ is a submanifold with dimension $\dim E$, and the slices at $a$ and $b$ are the parameter-boundary faces of this finite-dimensional reduction.

By the definition, on $M$ the map $F$ takes values in $E$. Let $f=F|_{M}:M\to E$. The zero set of $f$ is $Z$. Since $M$ is finite-dimensional and $Z$ is compact, we can take a relatively compact neighborhood $D$ of $Z$ in $M$ such that $\min_{\partial_{{\rm lat}}D}\|f\|_{E} > 0$, where $D$ is relatively open in $M$ and $\partial_{{\rm lat}}D = \overline{D}^{M} \setminus D$. Moreover, $\min_{D_b}\|f\|_E > 0$ where $D_b = \overline{D} \cap M_b$. Using the assumption that $F(a,\cdot)$ has exactly one zero $x_a$, and $D_xF(a,x_a)$ is an isomorphism, it is straightforward to check that $f|_{M_a}$ has only one zero, which is also regular.

The ordinary finite-dimensional mod-$2$ degree now gives a contradiction. By a transverse perturbation, we may assume that $f$ is smooth on $\overline{D}$. In this case, the zero set $Z$ is a compact one-dimensional manifold. Since $Z \cap M_b = \emptyset$, we conclude that $Z$ has exactly one boundary point, a contradiction.
\end{proof}

\medskip

\noindent \textbf{Acknowledgements:} This work is funded by National Key R\&D Program of China (Grant 2023YFA1010001), NSFC(12571123) and NSFC(12671139). It was completed when Song was visiting the Universit\'e Marie et Louis Pasteur, LmB (UMR 6623). Song is supported by the China Postdoctoral Science Foundation (2024T170452).

\end{document}